\documentclass[10pt]{article}

\usepackage{amssymb,amsmath,amsthm}
\usepackage{amsfonts}
\usepackage{graphicx}
\usepackage[usenames]{color}
\usepackage{caption}
\usepackage{bm,mathrsfs}
\usepackage{url}
\usepackage{subfigure}
\usepackage{algorithm,algorithmic}
\usepackage{dsfont,verbatim}
\usepackage{epsfig}
\usepackage[margin=1.2in]{geometry}
\usepackage{enumerate}

\graphicspath{{./pics/}}
\allowdisplaybreaks
\newtheorem{theorem}{Theorem}[section]
\newtheorem{lemma}{Lemma}[section]

\newtheorem{remark}{Remark}[section]

\newtheorem{proposition}{Proposition}[section]
\numberwithin{equation}{section}

\def\d{{\rm d}}

\title{\bf  Error analysis of energy-conservative BDF2--FE scheme for the 2D
Navier--Stokes equations with variable density}

\author{Jingjing Pan
\thanks{Department of Mathematics, School of Sciences, Hangzhou Dianzi University, Hangzhou, Zhejiang, P. R. China.Email address: \texttt{211070032@hdu.edu.cn}}
\and Wentao Cai
\thanks{Department of Mathematics, School of Sciences, Hangzhou Dianzi University, Hangzhou, Zhejiang, P. R. China.Email address: \texttt{femwentao@hdu.edu.cn}}
}

\date{}

\begin{document}

\maketitle

\begin{abstract}
In this paper, we present an error estimate of a second-order linearized finite element (FE) method for the 2D Navier--Stokes equations with variable density. In order to get error estimates, we first introduce an equivalent form of the original system. Later, we propose a general BDF2-FE method for solving this equivalent form,  where the Taylor--Hood FE space is used for discretizing the Navier--Stokes equations and conforming FE space is used for  discretizing density equation. We show that our scheme ensures discrete energy dissipation. Under the assumption of sufficient smoothness of strong solutions, an error estimate is presented for our numerical scheme for variable density incompressible flow in two dimensions. Finally, some numerical examples are provided to confirm our theoretical results.
\\[5pt]
\end{abstract}

\pagestyle{myheadings}
\thispagestyle{plain}

\section{Introduction}\label{sec:intro}

We consider the following Navier--Stokes equations with variable density in a bounded convex polygon domain $\varOmega\subset\mathbb{R}^2$, for $t\in(0,T)$
\begin{align}
\partial _t\rho + \nabla  \cdot (\rho {\bf u}) =& 0, \label{ii1} \\
\rho({{\bf u}_t} +{\bf u} \cdot\nabla {\bf u})-\mu\Delta {\bf u} + \nabla p =& 0, \label{ii2}\\
\nabla  \cdot{\bf u}  =& 0.\label{ii3}
\end{align}
In this system, ${\bf u}$ and $p$ represent the velocity and pressure of the fluid,
respectively. $\rho$ denotes the density.
Here we assume the system \eqref{ii1}-\eqref{ii3} subjects to the following initial
and boundary conditions:
\begin{equation}\label{ii4}
\left\{
     \begin{aligned}
&{\bf u} = {\bf 0},&&\,\,\, \mbox{on}\hspace{0.1cm} \partial\varOmega\times[0,T],\\
&\rho\,(x,0) = \rho^0\hspace{0.1cm}\mbox{and}\hspace{0.1cm} {\bf u}\,(x,0) = {\bf u}^0,&&\,\,\hspace{0.1cm}\mbox{in} \hspace{0.1cm}\varOmega,
     \end{aligned}
\right.
\end{equation}
where $\rho^0$ and ${\bf u}^0$ are two given functions, and $\partial\varOmega$ is denoted by the boundary of domain. We assume that initial date $\rho^0$ and ${\bf u}^0$ are sufficiently smooth, and
\begin{align}\label{minrho}
& \min \limits_{x \in \varOmega } \rho^0(x)>0.
\end{align}

The incompressible Navier--Stokes equations with variable density have been widely used in the fluid dynamics, such as highly stratified flows, interfaces between fluids of different densities, and astrophysics. The mathematical theory of Navier--Stokes equations with variable density were studied in \cite{RD,OVU,EPSM} for the existence and uniqueness of smooth solutions in two dimensions.

Numerical methods for the variable density Navier--Stokes equations have been investigated extensively in the past several decades, such as the projection  method \cite{ASA,JBB,JLL,YLLJF},
fractional-step method \cite{AAA,BBB}, backward difference method \cite{YLI},
discontinuous Galerkin (DG) method \cite{CLNJ,DGli}, and Gauge--Uzawa method in \cite{GU1,JShen}.
In \cite{JBB}, Bell and Marcus proposed a second-order projection method for the variable density Navier--Stokes equations, where the Crank--Nicolson method was used for discretization of time, and a standard difference method was used for discretization of space. Some numerical results are presented to varify the convergence properties.
A new time-stepping method was proposed in \cite{BBB}. Different from other classic
algorithms, the discrete pressure variable was solved by introducing one Possion
equation, which leads to reduce the computational cost.
Li et al. in \cite{YLI} proposed a second-order mixed stabilized FE method to solve variable density Navier--Stokes equations, where a second-order
backward difference was used for time discretization and the nonlinear term was explicitly treated by a second-order Adams--Bashforth method. The stability of
this scheme was proved, but no error estimates was presented.
Liu and Walkington \cite{CLNJ} investigated the DG method for the variable density
Navier--Stokes equations. The convergence of scheme was proved and no convergence
rates was presented.
Pyo and Shen in \cite{JShen} studied two Gauge--Uzawa schemes. The authors showed that the first-order temporal discretized Gauge--Uzawa schemes are unconditionally
stable.

Due to the strong nonlinearities and coupling terms,
it is rather difficult for error analysis of numerical methods for Navier--Stokes
equations with variable density. Recently, Guermond and Salgado in  \cite{JJLG}  proposed a fractional time-stepping method to solve Navier--Stokes equations with variable density. They provided the error estimate for a velocity equation \eqref{ii1} by the assumption of uniform boundedness of numerical density variable
 $\rho_h^{n}$, $n=1, \ldots,N.$   Cai et al. \cite{WBY} presented the first complete error estimate of the backward Euler method for the 2D Navier--Stokes equations with variable density by a error splitting technique and discrete maximal $L^p$-regularity. However, the analysis in \cite{WBY} could not  extended to the 3D problem. In \cite{DGli}, Li, Qiu and Yang proposed a linearized semi-implicit finite element
method, which utilized $H^1$-conforming finite element method to solve velocity equation and DG finite element method with post-processed numerical velocity variable to solve density equation.
The authors provided a complete error estimate of this numerical scheme in a three-dimensional domain.
Based on this work, Li and An in \cite{LiBDF2} proposed a BDF2-FE scheme, where the MINI element space was used for approximating velocity
${\bf u}$ and pressure $p$, and the quadratic conforming finite element space was
used for approximating density $\rho$. By a post-processed technique \cite{DGli},
the optimal convergence order $O(\tau^2+h^2)$ for numerical density variable
$\rho_h$ and numerical velocity variable ${\bf u}_h$ was proved in $L^2$-norm.
But, there is no literature on error estimates of the general (without post-processing for velocity) BDF2-FE scheme for Navier--Stokes equations with variable
density.

In this paper, we will present a error estimate of BDF2-FE scheme.
First, in order to analyze, we rewrite the system \eqref{ii1}-\eqref{ii3} as the
following equivalent system:
\begin{align}
\partial_t \sigma + {\bf u} \cdot \nabla \sigma +\frac{1}{2} \sigma \nabla \cdot {\bf u}=&0,
\label{eq1}\\
\sigma(\sigma {\bf u})_t +\rho ({\bf u}\cdot \nabla {\bf u})
+\frac{{\bf u}}{2}\nabla\cdot (\rho {\bf u}) -\mu\Delta {\bf u} + \nabla p=&0,
\label{eq2}\\
\nabla  \cdot {\bf u}  = &0,
\label{eq3}
\end{align}
where $\sigma =\sqrt{\rho}$.
Then we investigate a general BDF2 finite element method for the above system, where we combine the Taylor--Hood (P2--P1) finite element space for discretizing velocity ${\bf u}$ and pressure $p$, and conforming finite element space (P2) for density
$\rho$. We present that our scheme keeps energy dissipation. Under the assumption
of sufficiently smooth of strong solution $(\rho, {\bf u},p)$, we present a complete
error estimate for a second-order temporal discrete FE scheme for Navier--Stokes equations with variable density. To our best knowledge, this paper is the first work
for a complete error estimate of the general BDF2--FEM for Navier--Stokes equations with variable density.

The rest of this paper is organized as follows. In Sect.2, we introduce the notations, numerical scheme and   the main results .  The proof of the main theorem is presented in Sect.3.
At last, numerical examples are presented in Sect.4 to support the theoretical analysis.

\section{Main results}
\subsection{Notations and variational formulation}\label{Sect:notation}
For any positive integer $k$ and real number $p$, we denote by $W^{k,p}(\varOmega)$ the conventional Sobolev space of functions defined on the domain $\varOmega$, with the abbreviations
 $H^k(\varOmega) = W^{k,2}(\varOmega)$, $L^p(\varOmega)=W^{0,p}(\varOmega)$ and
$L_0^p(\varOmega)=\{q\in L^{p}(\varOmega):\mbox{$\int_{\varOmega}qdx=0$}\}$. The space of continuous functions on $\overline\varOmega$ is denoted by $C(\overline\varOmega)$.
The closure of $C_0^{\infty}(\varOmega)$ in $W^{k,p}(\varOmega)$ space is denoted
by $W_0^{k,p}(\varOmega)$,
with the abbreviation $H^{k}_0(\varOmega)= W^{k,2}_0(\varOmega)$.
The boldface notation $\textbf{W}^{k,p}(\varOmega)$
and $\textbf{L}^p(\varOmega)$ are used to denote the vector-value Sobolev spaces corresponding to  $[W^{k,p}(\varOmega)]^2$ and $[L^p(\varOmega)]^2$,
with the abbreviations ${\bf H}^{k}(\varOmega)={\bf W}^{k,2}(\varOmega)$ and $ {\bf L}^p(\varOmega)= {\bf W}^{0,p}(\varOmega)$.

For simplicity, we denote the inner products of both
$L^2(\varOmega)$ and $\textbf{L}^2(\varOmega)$
by $(\cdot,\cdot)$, namely,
\begin{align*}
&(u,v)=\int_\varOmega u(x)v(x)\d x, &&\forall\,\, u,v\in L^2(\varOmega),\\
&({\bf u},{\bf v})=\int_\varOmega {\bf u}(x)\cdot {\bf v}(x)\d x, &&\forall\, \,{\bf u},{\bf v}\in {\bf L}^2(\varOmega) .
\end{align*}

With above notations, it is easily seen that the exact solution
$(\sigma, {\bf u}, p)$
of \eqref{eq1}-\eqref{eq3} satisfies the following variation equations:
\begin{align}
&
(\partial_t \sigma,\varphi) + (\nabla \sigma\cdot {\bf u},\varphi) +\frac{1}{2} (\sigma \nabla \cdot {\bf u},\varphi)=0,
\label{eq4}\\
&
(\sigma(\sigma {\bf u})_t,{\bf v}) +(\rho ({\bf u} \cdot \nabla {\bf u}),{\bf v})+\frac{1}{2}(\nabla\cdot (\rho {\bf u})\cdot {\bf u},{\bf v})
+ B(({\bf u},p),({\bf v},q))=0,
\label{eq5}
\end{align}
for any test functions $\varphi\in L^2(\varOmega)$ and $({\bf{v}},q)\in {\bf{H}}_0^1(\varOmega)\times L_0^2(\varOmega)$, where $B(({\bf u},p),({\bf v},q))$
is the bilinear form defined by:
\begin{align}
B(({\bf u},p),({\bf v},q))=(\mu\nabla {\bf u},\nabla {\bf v})
-(p,\nabla\cdot {\bf v})+(\nabla\cdot {\bf u},q),
\hspace{0.4cm} \forall\,\,({\bf u},p),({\bf v},q) \in {\bf H}_0^1(\varOmega)\times L_0^2(\varOmega).
\nonumber
\end{align}

In the equation above, the extra term
$\frac{1}{2} (\sigma \nabla \cdot {\bf u},\varphi)$ vanish due to \eqref{eq3}. This extra term will stabilize the finite element solutions to be defined in Subsection 2.3.

\subsection{Numerical scheme and main results}

Let $\{t_n=n\tau \}_{n=0}^N$ be a uniform partition of the time interval $[0,T]$ with time step $\tau=T/N$. For any sequence $\{z^n\}_{n=0}^N $,
we denote
\begin{align}
&\hat{z}^{n+1}=2z^n-z^{n-1},
\hspace{0.8cm}
D_\tau z^{n+1}=\frac{ 3z^{n+1}-4z^n+z^{n-1}}{2\tau},
\,\,\,\,\, n=1,2,\ldots,N-1.
\nonumber \\
&\delta_\tau z^{n+1}=\frac{z^{n+1}-z^n}{\tau},\,\,\,n=0,1,\ldots,N-1.
\nonumber
\end{align}
Let $\mathscr{T}_h$ be a quasi-uniform triangulation of $\varOmega$ into triangles ${\mathcal T}_j$, $j=1,...,M$, with mesh size $h=\max_{{1\leq j\leq M}}{\rm diam}({\mathcal T}_j)$.
For any integer $r\geq1$, we define the following finite element spaces:
\begin{align}
&M_h^r=\{u_h\in H^1(\varOmega): u_h|_{{\mathcal T}_j}\in P_r({\mathcal T}_j), \hspace{0.3cm} \forall\,\,\, {\mathcal T}_j\in \mathscr{T}_h\},
\nonumber \\
&{\bf X}_h^{r}=\{ p_h\in H_0^1(\varOmega)^2:
p_h|_{{\mathcal T}_j}\in P_r({\mathcal T}_j)^2, \hspace{0.1cm} \forall\,\,\,
{\mathcal T}_j\in \mathscr{T}_h \},
\nonumber \\
&\mathring M_h^1=\{q_h\in M_h^1:  \mbox{$\int_\varOmega q_h(x)\d x =0$}\},
\nonumber
\end{align}
where $P_r({\mathcal T}_j)$ is the space of  polynomials of degree $\le r$ on the triangle ${\mathcal T}_j$.
Consequently, ${\bf X}_h^{2}\times \mathring M_h^1$ is a Taylor--Hood finite element space,
which satisfies the following inf-sup condition (cf.\cite{DBOff})
\begin{align}
\|q_h\|_{L^2(\varOmega)}\le C \sup_{{\bf v}_h \in {\bf X}_h^{2}} \frac{|(\nabla q_h,{\bf v}_h)|}{\|\nabla {\bf v}_h\|_{L^2(\varOmega)}},
\quad
\forall\, q_h\in \mathring M_h^1 .
\label{infsupcond}
\end{align}

In this paper, we use $C$ to represent generic constants, while $\epsilon$ denotes a small generic constant. Both $C$ and $\epsilon$ are independent of $\tau$, $h$ and $N$, and their values may differ in various instances.

Based on the above notations, a stabilized BDF2 finite element scheme for the
Navier--Stokes equations with variable density is to find
$(\rho_h^{n+1}, {\bf u}_h^{n+1}, p_h^{n+1})\in M_h^2\times \textbf{X}_h^{2} \times \mathring M_h^1$, such that
\begin{align}
&(D_\tau\sigma_h^{n+1},\varphi_h)
+(\nabla\sigma_h^{n+1}\cdot\hat{{\bf u}}_h^{n+1},\varphi_h)
+\frac{1}{2}(\nabla\cdot\hat{{\bf u}}_h^{n+1}\sigma_h^{n+1},\varphi_h)=0,
\label{li1}\\
&(\sigma_h^{n+1}D_\tau(\sigma_h^{n+1}{\bf u}_h^{n+1}),{\bf v}_h)
+(\rho_h^{n+1}\hat{{\bf u}}_h^{n+1}\cdot\nabla{\bf u}_h^{n+1},{\bf v}_h)
\nonumber \\
&\,\,\,
+\frac{1}{2}({\bf u}_h^{n+1}\nabla\cdot(\rho_h^{n+1}\hat{{\bf u}}_h^{n+1}),{\bf v}_h)
+B(({\bf u}_h^{n+1},p_h^{n+1}),({\bf v}_h,q_h))=0,
\label{li2}
\end{align}
hold for all test functions
$(\varphi_h,{\bf v}_h,q_h)\in M_h^2\times \textbf{X}_h^{2} \times \mathring M_h^1$, where $\sigma_h^{n+1}=\sqrt{ \rho_h^{n+1}}$.

In this paper, the following inverse inequality is frequently used:
\begin{align}
\|v_h\|_{W^{m,s}}\leq Ch^{n-m+\frac2s-\frac2q}\|v_h\|_{W^{n,q}},
\label{inverse}
\end{align}
for any $v_h\in V_h$, $0\le n\le m\le 1$ and $1\le q\le s\le \infty$ in two dimensions.

\begin{remark}\label{inital}
The starting values
$(\rho_h^i,{\bf u}_h^i,p_h^i)\in M_h^2\times \textbf{X}_h^{2} \times \mathring M_h^1,$ $i=0,1$ are assumed to be given and satisfy the estimates \eqref{F1}-\eqref{F3}. An example of constructing the numerical schemes for starting values is given in Section 4, which
satisfies the convergence assumed in \eqref{F1}-\eqref{F3}.
\end{remark}

In order to get error estimate, we need the following assumption
conditions:
\begin{align}
  &\mbox{The solution of \eqref{ii1}-\eqref{ii4} is sufficiently smooth.
  \qquad\qquad\qquad\qquad\qquad\qquad\qquad\qquad\quad\,\,\,\,}
  \label{smooth}  \\
  &\mbox{The starting values $(\rho_h^i,{\bf u}_h^i,p_h^i)\in M_h^2\times  \textbf{X}_h^{2} \times \mathring M_h^1,$ $i=0,1$ are given sufficiently accurate, }
  \nonumber \\
  &\mbox{ satisfying the following estimates:}
  \nonumber \\
  &\qquad\qquad\qquad\qquad\qquad\qquad\qquad\quad\|P_h\sigma^i-\sigma_h^i\|_{L^2}\leq C(\tau^3+h^2) ,
  \label{F1} \\
  &\qquad\qquad\qquad\qquad\qquad\qquad\qquad\quad \|{\bf R}_h{\bf u}^i-{\bf u}_h^i\|_{L^2}
  \leq C(\tau^2+h^2),
  \label{F2}\\
  &\qquad\qquad\qquad\qquad\qquad\qquad\qquad\quad\tau\|{\bf R}_h{\bf u}^i-{\bf u}_h^i\|_{H^1}^2
  \leq Ch^2(\tau^2+h^2),
  \label{F3}
\end{align}
for $i=0,1$.
In \eqref{F1}-\eqref{F3}, ${\bf R}_h$ and $P_h$ denote the Ritz and Stokes projection, respectively. The definitions of these projections can be found in Subsection 3.1.

In this paper, we assume that the initial density $\rho^0 > 0$  satisfies
\begin{align}
&
0\leq \min_{x\in\overline\varOmega}\rho^0(x)\leq\max_{x\in\overline\varOmega}\rho^0(x)
< \infty,
\nonumber
\end{align}
and recall that (cf. \cite{AN} and \cite{rho})
\begin{align}
\mbox{meas} \,\,\{ x\in\varOmega, \alpha\leq\rho(x,t)\leq\beta\}
\,\,\mbox{is independent of} \,\,t \geq 0
\,\,\mbox{for all}\,\, 0\leq\alpha\leq\beta<+\infty,
\nonumber
\end{align}
then, we get
\begin{align}
\min_{x\in\overline\varOmega}\rho^0(x)\leq\rho(x,t)\leq\max_{x\in\overline\varOmega}\rho^0(x),\qquad\,\,
\forall \,\,t\geq 0,
\nonumber
\end{align}
which implies that
\begin{align}
\min_{x\in\overline\varOmega}\sigma^0(x)
\leq\sigma(x,t)\leq\max_{x\in\overline\varOmega}\sigma^0(x)
,\qquad\,\,
\forall \,\,t\geq 0.
\label{prove-sig}
\end{align}

Moreover, we will use the following discrete Gronwall's
inequality in error analysis:
\begin{lemma}\label{Gronwall's}
(Discrete Gronwall's inequality) let $\tau$, B, and $a_k$, $b_k$, $c_k$, $\gamma_k$,
 for integers $k\geq0$, be non-negative numbers such that
 \begin{align}\
 \begin{aligned} \nonumber
 a_j+\tau\sum\limits_{k=0}^j b_k \leq
 \tau\sum\limits_{k=0}^j\gamma_k a_k
 + \tau \sum\limits_{k=0}^j c_k + B,
 &&
 \textrm{for} \hspace{0.2cm} j\geq0,
 \end{aligned}
 \end{align}
suppose that $\tau\gamma_k<1$, for all k, and set
$\vartheta_k=(1-\tau\gamma_k)^{-1}$. Then
\begin{align}
\begin{aligned}\nonumber
a_j+\tau\sum\limits_{k=0}^j b_k
\leq \exp(\tau\sum\limits_{k=0}^j \gamma_k\vartheta_k)
(\tau\sum\limits_{k=0}^j c_k+B),
&&
\textrm{for}\hspace{0.2cm}j\geq0.
\end{aligned}
\end{align}
\end{lemma}

The main theoretical result of this paper is presented as follows:

\begin{theorem}\label{conclusion}
Assume that the system \eqref{ii1}-\eqref{minrho} has a unique solution
$(\rho,{\bf u},p)$ that satisfies the regularity assumption \eqref{smooth}-\eqref{F3},
then the fully discrete system \eqref{li1}-\eqref{li2} yields a unique solution
$(\rho_h^{n+1},{\bf u}_h^{n+1}, p_h^{n+1})$. Moreover, there exist positive constants
$\tau_*$ and $h_*$ such that when $\tau\leq\tau_*$,  $h\leq h_*$ and $\tau = O(h)$,
 satisfying
\begin{align}\label{theorem}
\max\limits_{1\leq n\leq N-1}&\left(
\|{\bf u}^{n+1}-{\bf u}_h^{n+1}\|_{L^2}+\|\rho^{n+1}-\rho_h^{n+1}\|_{L^2}\right)
\leq C(\tau^2 + h^2).
\end{align}
\end{theorem}

\subsection{Unconditional energy stability}

In this section, we present the unconditional energy stability of our numerical scheme \eqref{li1}-\eqref{li2}. Further, we will design numerical examples to
verify this property in Section 4.

\begin{proposition}\label{energy}
The FE scheme \eqref{li1}-\eqref{li2} ensures the energy dissipation, that is
\begin{align}
E^{n+1}\leq E^n,  \,\,\,\,\, n=1,2,\ldots,N-1,
\nonumber
\end{align}
where the discrete energy $E^n$ is defined by
\begin{align}
E^{n}\buildrel \Delta \over =
\|\sigma_h^{n}\|_{L^2}^2
+\|\sigma_h^{n}{\bf u}_h^{n}\|_{L^2}^2
+\|2\sigma_h^{n}-\sigma_h^{n-1}\|_{L^2}^2
+\|2\sigma_h^{n}{\bf u}_h^{n}-\sigma_h^{n-1}{\bf u}_h^{n-1}\|_{L^2}^2.
\label{EEE}
\end{align}
\end{proposition}

\noindent
{\bf Proof}\,\,\,
Taking $\varphi_h=\sigma_h^{n+1}$ and
$({\bf v}_h,q_h)=( {\bf u}_h^{n+1},p_h^{n+1})$ in \eqref{li1} and \eqref{li2}
respectively, we obtain
\begin{align}
&(D_\tau\sigma_h^{n+1},\sigma_h^{n+1})
+(\nabla\sigma_h^{n+1}\cdot\hat{{\bf u}}_h^{n+1},\sigma_h^{n+1})
+\frac{1}{2}(\nabla\cdot\hat{{\bf u}}_h^{n+1}\sigma_h^{n+1},\sigma_h^{n+1})
=0,
\label{EN1}\\
&(\sigma_h^{n+1}D_\tau(\sigma_h^{n+1}{\bf u}_h^{n+1}),{\bf u}_h^{n+1})
+(\rho_h^{n+1}\hat{{\bf u}}_h^{n+1}\cdot\nabla{\bf u}_h^{n+1},{\bf u}_h^{n+1})
\nonumber \\
&\,\,\,
+\frac{1}{2}({\bf u}_h^{n+1}\nabla\cdot(\rho_h^{n+1}\hat{{\bf u}}_h^{n+1}),{\bf u}_h^{n+1})
+B(({\bf u}_h^{n+1},p_h^{n+1}),({\bf u}_h^{n+1},p_h^{n+1}))
=0.
\label{EN2}
\end{align}
According to  the integration by parts, we get
\begin{align}
(\nabla\sigma_h^{n+1}\cdot\hat{{\bf u}}_h^{n+1},\sigma_h^{n+1})
&=-\frac{1}{2}(\nabla\cdot\hat{{\bf u}}_h^{n+1}\sigma_h^{n+1},\sigma_h^{n+1}),
\nonumber \\
(\rho_h^{n+1}\hat{{\bf u}}_h^{n+1}\cdot\nabla{\bf u}_h^{n+1},{\bf u}_h^{n+1})
&=
-\frac{1}{2}({\bf u}_h^{n+1}\nabla\cdot(\rho_h^{n+1}\hat{{\bf u}}_h^{n+1}),{\bf u}_h^{n+1}),
\nonumber
\end{align}
which with \eqref{EN1}, \eqref{EN2} and the fact
\begin{align}\label{KKK}
2(a^{k+1}, 3a^{k+1}-4a^k+a^{k-1})=&|a^{k+1}|^2-|a^{k}|^2+|a^{k+1}-2a^k+a^{k-1}|^2
\nonumber \\
&+|2a^{k+1}-a^{k}|^2-|2a^{k}-a^{k-1}|^2,
\end{align}
imply
\begin{align}
\frac{1}{4\tau}&
\left(\|\sigma_h^{n+1}\|_{L^2}^2-\|\sigma_h^{n}\|_{L^2}^2
+\|2\sigma_h^{n+1}-\sigma_h^{n}\|_{L^2}^2
-\|2\sigma_h^{n}-\sigma_h^{n-1}\|_{L^2}^2
\right)
\leq0,
\\
\frac{1}{4\tau}
&\Big(\|\sigma_h^{n+1}{\bf u}_h^{n+1}\|_{L^2}^2-\|\sigma_h^{n}{\bf u}_h^{n}\|_{L^2}^2
+\|2\sigma_h^{n+1}{\bf u}_h^{n+1}-\sigma_h^{n}{\bf u}_h^{n}\|_{L^2}^2
\nonumber \\
&\qquad\qquad\quad\,\,\,\,\,\,\,\,\,\,\,\,
-\|2\sigma_h^{n}{\bf u}_h^{n}-\sigma_h^{n-1}{\bf u}_h^{n-1}\|_{L^2}^2\Big)
+\mu\|\nabla {\bf u}_h^{n+1}\|_{L^2}^2\leq0.
\end{align}
Thus, we get $E^{n+1}\leq E^n$, $n = 1,\ldots, N-1$. The proof of
Proposition \ref{energy} is complete.

\section{Error analysis of the fully discrete scheme \eqref{li1}-\eqref{li2}}

\subsection{Ritz projection, ${\bf L^2}$ projections, and Lagrange interpolations }

Let $(\textbf{R}_h,Q_h):{\bf{H}}_0^1(\varOmega)\times L^2(\varOmega)\rightarrow \textbf{X}_h^{2}\times \mathring M_h^1$ be a class Stokes--Ritz projection defined by
\begin{equation}
B \left((\textbf{R}_h({\bf u},p),Q_h({\bf u},p)),({\bf v}_h,q_h)\right)=B(({\bf u},p),({\bf v}_h,q_h)),
\hspace{0.2cm}\forall\,\,({\bf v}_h,q_h)\in {\bf{X}}_h^{2}\times M_h^1.
\nonumber
\end{equation}
It is well known that the Ritz projection satisfies the following standard estimates
(cf. \cite{Ritz}, Thm. 12.6.7, \cite{Ritz-L2}, Chap. II, Thm. 4.3, \cite{Ritz-W1}, Thm. 8.2):
\begin{align}
&
\|\textbf{R}_h {\bf v}-{\bf v}\|_{L^2}\leq Ch^{l+1}(\|{\bf v}\|_{H^{l+1}}+\|q\|_{H^l}),
\qquad\qquad\qquad\qquad\quad\,\,\,
 l=1,2,
\label{T1}\\
&
\|\textbf{R}_h {\bf v}-{\bf v}\|_{H^1} +\|Q_h q-q\|_{L^2} \leq
Ch^l(\|{\bf v}\|_{H^{l+1}}+\|q\|_{H^{l}}),
\qquad\quad\,\,\,\,
 l=1,2,
\label{T2} \\
&
\|\textbf{R}_h {\bf v}\|_{W^{1,\infty}} \leq
C(\|{\bf v}\|_{W^{1,\infty}}+\|q\|_{L^\infty}).
\label{T3}
\end{align}
for
$({\bf v},q)\in( \textbf{H}^{l+1} ( \varOmega)\cap {\bf{H}}_0^1(\varOmega))\times H^l(\varOmega)$. In the above estimates, we denote $\textbf{R}_h {\bf v}:=\textbf{R}_h( {\bf v},q)$,
$Q_h q:=Q_h({\bf v},q)$ for simplicity.

Let $P_h : L^2(\varOmega)\rightarrow M_h^2$ and ${\bf P}_h : {\bf L}^2(\varOmega)\rightarrow {\bf X}_h^{2}$ denote the $L^2$ projection operators defined by
\begin{align}
&(v-P_h v, w_h)=0,
&&\forall \,\,w_h\in M_h^2,
\label{VV1}\\
&({\bf v}-P_h {\bf v}, {\bf w}_h)=0,
&&\forall \,\,{\bf w}_h\in {\bf X}_h^{2}.
\label{VV2}
\end{align}
By classic FE theory \cite{Ritz}, the $L^2$ projection satisfies the following estimates:
\begin{align}
&
\|\varphi-P_h\varphi\|_{L^2} + h\|\nabla(\varphi-P_h\varphi)\|_{L^2}
 \leq Ch^{l+1}\|\varphi\|_{H^{l+1}},
&&
l=1,2,
\label{T4}\\
&
\|P_h \varphi\|_{W^{k,q}}\leq C \|\varphi\|_{W^{k,q}},
&&k=0,1, \hspace{0.2cm}1\leq q \leq \infty.
\label{T6}
\end{align}

Similarly, the Lagrangian interpolation operators
$\Pi_h : C(\varOmega)\rightarrow M_h^2$ and
${\bf \Pi}_h : {\bf C}(\varOmega)\rightarrow {\bf X}_h^{2}$ satisfy:
\begin{align}
&
\|\Pi_h\varphi-\varphi\|_{L^2} + h \|\nabla(\Pi_h\varphi-\varphi)\|_{L^2}\leq Ch^{l+1}\|\varphi\|_{H^{l+1}},
&&\forall\,\, \varphi \in  H^{l+1}(\varOmega),
&&l=1,2,
\label{pi0}\\
&
\|{\bf \Pi}_h{\bf{v}}-{\bf{v}}\|_{L^2} + h \|\nabla({\bf \Pi}_h{\bf{v}}-{\bf{v}})\|_{L^2}\leq Ch^{l+1}\|{\bf{v}}\|_{H^{l+1}},
&&\forall\,\, {\bf{v}} \in  {\bf{H}}^{2}(\varOmega) \cap {\bf{H}}_0^1(\varOmega),
&&l=1,2,
\label{pi1}\\
&
\|{\bf \Pi}_h{\bf{v}}-{\bf{v}}\|_{L^\infty} \leq Ch \|{\bf{v}}\|_{W^{1,\infty}},
&&\forall\,\, {\bf{v}} \in  {\bf{W}}^{1,\infty}(\varOmega) \cap {\bf{H}}_0^1(\varOmega).
\label{pi2}
\end{align}
The estimates \eqref{T1}-\eqref{T3} and \eqref{T4}-\eqref{pi2} will be frequently used in this section.

\subsection{The proof of Theorem 2.1}

In the following part, we prove the existence and uniqueness of the numerical solutions of
system \eqref{li1}-\eqref{li2}.
If $(\rho_h^k,{\bf u}_h^k, p_h^k)\in M_h^2\times \textbf{X}_h^{2} \times \mathring M_h^1$ is given for $k=0,1,2,\ldots, n$, then the FE equation \eqref{li1} has a unique solution if and only if the corresponding homogeneous equation
\begin{align}
\frac 3 2\tau^{-1}(\varTheta,\varphi_h)
+(\hat{{\bf u}}_h^{n+1}\cdot\nabla\varTheta,\varphi_h )
+\frac 1 2(\nabla\cdot\hat{{\bf u}}_h^{n+1}\varTheta,\varphi_h)=0,
\qquad\forall\,\,\varphi_h\in M_h^2,
\end{align}
has only zero solution $\varTheta=0$.
According to the integration by parts, the above equation  can be rewritten
as
\begin{align}
\frac 3 2\tau^{-1}(\varTheta,\varphi_h)
+\frac 1 2(\hat{{\bf u}}_h^{n+1}\cdot\nabla\varTheta,\varphi_h )
-\frac 1 2(\hat{{\bf u}}_h^{n+1}\varTheta,\nabla\varphi_h)=0.
\nonumber
\end{align}
Taking $\varphi_h=\varTheta$ in the above equation, we have $\|\varTheta\|_{L^2}=0$.
It implies $\varTheta=0$. Thus, we prove the unique solvability of equation
\eqref{li1}.

After we solve $\sigma_h^{n+1}(\sigma_h^{n+1}=\sqrt{\rho_h^{n+1}})$ in \eqref{li1}, the FE equation \eqref{li2} has a unique solution $({\bf u}_h^{n+1}, p_h^{n+1})$ if and only if the homogeneous equation
\begin{align}
&\frac 3 {2}\tau^{-1}(\sigma_h^{n+1}\sigma_h^{n+1}{\bf U}_h,{\bf v}_h)
+(\rho_h^{n+1}\hat{{\bf u}}_h^{n+1}\cdot\nabla {\bf U}_h,{\bf v}_h)
+\frac 1 2({\bf U}_h\nabla\cdot(\rho_h^{n+1}\hat{{\bf u}}_h^{n+1}),{\bf v}_h)
\nonumber \\wo
&\quad
+B(({\bf U}_h,P_h),({\bf v}_h,q_h))=0,
\label{UUUU}
\end{align}
has only zero solution $({\bf U}_h,P_h)=({\bf 0},0).$ Indeed, substituting
$({\bf v}_h,q_h)=({\bf U}_h,P_h)$ into the equation above yields
\begin{align}
\frac 3 2 \tau^{-1}\|\sigma_h^{n+1}{\bf U}_h\|_{L^2}^2
+\|\nabla{\bf U}_h\|_{L^2}^2=0,
\nonumber
\end{align}
which with the fact $\rho_h^{n+1}>0$, implies ${\bf U}_h=0$. Then \eqref{UUUU}
reduces to
\begin{align}
(P_h,\nabla\cdot {\bf v}_h)=0,\qquad\quad\forall\,\,{\bf v}_h\in{\bf X}_h^{2},
\nonumber
\end{align}
which implies $\|P_h\|_{L^2}=0$ by the inf-sup condition \eqref{infsupcond}.
Thus, we get the unique solvability of equation \eqref{li2}.

Next, we begin to prove the error estimate \eqref{theorem}.
For the numerical solution $(\rho_h^{n},{\bf u}_h^{n},p_h^{n})$ of \eqref{li1}-\eqref{li2}, we denote
\begin{align}\label{eee}
&
e_{\rho,h}^{n}=P_h\rho^{n}-\rho_h^{n},
&&e_{{\bf u},h}^{n}={\bf{R}}_h{\bf u}^{n}-{\bf u}_h^{n}, \nonumber
\\
&e_{\sigma,h}^{n}=P_h\sigma^{n}-\sigma_h^{n},
&&e_{p,h}^{n}=Q_h p^{n}-p_h^{n}.
\end{align}

By mathematical induction,  we prove the following primary estimates: If the temporal step $\tau$ and spatial step $h$ are sufficiently small, then
\begin{subequations}\label{CC1}
\begin{align}
\max_{0\leq k\leq {n}}\|e_{\sigma,h}^{k}\|_{L^\infty}\leq&
\frac 1 4 \min_{x\in\overline\varOmega}\sigma^0(x),
\label{CC1a}
\\
\max_{0\leq k\leq {n}}\|e_{{\bf u},h}^{k}\|_{L^2}\leq&
h(\tau^2+h^2)^{\frac 1 4},
\label{CC1b}
\\
\max_{0\leq k\leq {n}}\|e_{{\bf u},h}^{k}\|_{L^\infty}\leq&1,
\label{CC1c}
\\
\tau\sum_{k=0}^{n}\|e_{{\bf u},h}^{k}\|_{H^1}^2
\leq&
h^2\sqrt{\tau^2+h^2},
\label{CC1d}
\end{align}
\end{subequations}
for $n=1,2,\ldots,N$.
By the estimates  \eqref{F1}-\eqref{F3},
we can derive that when $\tau$ and $h$ are sufficiently small,
and $\tau=O(h)$,
\begin{subequations}\label{0000}
\begin{align}
\|e_{\sigma,h}^i\|_{L^\infty}
\leq& Ch^{-1}\|e_{\sigma,h}^i\|_{L^2}
\leq Ch^{-1}(\tau^3+h^2)
\leq\frac 1 4 \min_{x\in\overline\varOmega}\sigma^0(x),
\label{0000a}
\\
\|e_{{\bf u},h}^i\|_{L^2}
\leq&
C(\tau^2+h^2)
\leq h(\tau^2+h^2)^{\frac 1 4},
\label{0000b}
\\
\|e_{{\bf u},h}^i\|_{L^\infty}\leq&
Ch^{-1}\|e_{{\bf u},h}^i\|_{L^2}
\leq Ch^{-1}(\tau^2+h^2)\leq 1,
\label{0000c}
\\
\tau\sum_{k=0}^{1}\|e_{{\bf u},h}^k\|_{H^1}^2
\leq&Ch^2(\tau^2+h^2)
\leq h^2\sqrt{\tau^2+h^2},
\label{0000d}
\end{align}
\end{subequations}
for $i=0,1.$ Thus, \eqref{CC1a}-\eqref{CC1d} hold for $n=1$.

Now, we assume that \eqref{CC1a}-\eqref{CC1d} are valid for $n\leq m$,
and satisfying the following inequalities:
\begin{subequations}\label{as2}
\begin{align}
\max_{0\leq k \leq m}\|e_{\sigma,h}^{k}\|_{L^\infty}\leq&
\frac 1 4 \min_{x\in\overline\varOmega}\sigma^0(x),
\label{as2a}
\\
\max_{0\leq k\leq m}\|e_{{\bf u},h}^{k}\|_{L^2}\leq&
h(\tau^2+h^2)^{\frac 1 4},\label{as2b}
\\
\max_{0\leq k\leq m}\|e_{{\bf u},h}^{k}\|_{L^\infty}\leq&1,\label{as2c}
\\
\tau\sum_{k=0}^{m}\|e_{{\bf u},h}^{k}\|_{H^1}^2\leq&
h^2\sqrt{\tau^2+h^2}.\label{as2d}
\end{align}
\end{subequations}
Note that the induction assumption \eqref{as2a} implies that
\begin{align}\label{sigma_n1}
\begin{aligned}
\max_{0\leq k\leq m}\|\sigma^k-\sigma_h^k\|_{L^\infty}
\leq&
\max_{0\leq k\leq m}\|\sigma^k-P_h\sigma^k\|_{L^\infty}+
\max_{0\leq k\leq m}\|e_{\sigma,h}^k\|_{L^\infty}
 \\
\leq&
 Ch^{-1}\max_{0\leq k\leq m}\|\sigma^k-P_h\sigma^k\|_{L^2}+
\max_{0\leq k\leq m}\|e_{\sigma,h}^k\|_{L^\infty}
 \\
\leq&
Ch\max_{0\leq k\leq m}\|\sigma^{k}\|_{H^{2}}
+\frac 1 4 \min_{x\in\overline\varOmega}\sigma^0(x)
 \\
\leq&
\frac 1 2 \min_{x\in\overline\varOmega}\sigma^0(x),
\quad\quad\,\,\,\mbox{(when $h$ is sufficiently small)}
\end{aligned}
\end{align}
which with \eqref{prove-sig} implies
\begin{align}
\begin{aligned}
&\sigma_h^{n}\geq\min_{x\in\overline\varOmega}\sigma^0(x)
-\|\sigma^{n}-\sigma_h^n\|_{L^\infty}
\geq\frac1 2 \min_{x\in\overline\varOmega}\sigma^0(x),
\quad\,\,\,\,\forall\,\, x\in\varOmega,\quad n=0,\ldots,m,
\nonumber \\
&\sigma_h^{n}\leq\max_{x\in\overline\varOmega}\sigma^0(x)
+\|\sigma^{n}-\sigma_h^n\|_{L^\infty}
\leq\frac3 2 \max_{x\in\overline\varOmega}\sigma^0(x),
\quad\,\,\forall \,\, x\in\varOmega,\quad n=0,\ldots,m,
\nonumber
\end{aligned}
\end{align}
thus
\begin{align}\label{sigma_n2}
\frac 1 2 \min_{x\in\overline\varOmega}\sigma^0(x)
\leq\sigma_h^n\leq\frac 3 2 \max_{x\in\overline\varOmega}\sigma^0(x),
\quad\quad n=0,\ldots,m.
\end{align}
Besides, the induction assumption \eqref{as2c} implies that
\begin{align}\label{U_n}
\|\hat{{\bf u}}&_h^{n+1}\|_{L^\infty}
\nonumber \\
&\leq\|{\bf{R}}_h \hat{{\bf u}}^{n+1} \|_{L^\infty}
+\|\hat{e}_{{\bf u},h}^{n+1}\|_{L^\infty}
\nonumber \\
&\leq\|{\bf{R}}_h \hat{{\bf u}}^{n+1} -\Pi_h\hat{{\bf u}}^{n+1}\|_{L^\infty}
+\|\Pi_h\hat{{\bf u}}^{n+1}-\hat{{\bf u}}^{n+1}\|_{L^\infty}+\|\hat{{\bf u}}^{n+1}\|_{L^\infty}
+\|\hat{e}_{{\bf u},h}^{n+1}\|_{L^\infty}
\nonumber  \\
&\leq
Ch^{-1}\|{\bf{R}}_h \hat{{\bf u}}^{n+1} -\Pi_h\hat{{\bf u}}^{n+1}\|_{L^2}
+\|\Pi_h\hat{{\bf u}}^{n+1}-\hat{{\bf u}}^{n+1}\|_{L^\infty}
+\|\hat{{\bf u}}^{n+1}\|_{L^\infty}
+\|\hat{e}_{{\bf u},h}^{n+1}\|_{L^\infty}
\nonumber \\
&\leq
Ch \left(\|\hat{{\bf u}}^{n+1}\|_{H^2}+\|p^{n+1}\|_{H^1}\right)
+Ch\|\hat{{\bf u}}^{n+1}\|_{W^{1,\infty}}+\|\hat{{\bf u}}^{n+1}\|_{L^\infty}+1
\nonumber \\
&\leq C,
\qquad\,\,\,\mbox{(when $h$ is sufficiently small)}
\end{align}
for $n=1,2,\ldots,m$.
Similar to the above induction, we can easily obtain
\begin{align}\label{U_n1}
\|{\bf u}_h^n\|_{L^\infty}\leq C,\qquad\quad
n=0,1,\ldots,m.
\end{align}

The estimates \eqref{sigma_n2}-\eqref{U_n1} will be frequently used in the following subsections.
Next, we want to prove that \eqref{CC1a}-\eqref{CC1d} also hold for $n=m+1$.

\subsection{Estimates of $e_{\sigma,h}^{n+1}$}
The equation \eqref{eq4} can be written as
\begin{align}
(D_\tau\sigma^{n+1},\varphi)+(\nabla \sigma^{n+1}\hat{{\bf u}}^{n+1},\varphi)
+\frac{1}{2}(\nabla \cdot \hat{{\bf u}}^{n+1}\sigma^{n+1},\varphi)
=(R_1^{n+1},\varphi),
\quad\quad\forall\,\,\varphi\in L^2(\varOmega),
\label{yuan1}
\end{align}
where $R_1^{n+1}$ is truncation error.  Under the regularity assumption \eqref{smooth}, we have
\begin{align}
\|R_1^{n+1}\|_{L^2}^2+\tau\sum_{n=0}^{N-1}\|R_1^{n+1}\|^2_{L^2}\leq C\tau^4.
\label{R1}
\end{align}
Subtracting \eqref{yuan1} from \eqref{li1} and
noting the $L^2$ projection \eqref{VV1} yields
\begin{align}\label{epp1}
(D_\tau &e_{\sigma,h}^{n+1} ,\varphi_h)
+
(\hat{{\bf u}}_h^{n+1}\nabla(\sigma^{n+1}-P_h\sigma ^{n+1}),\varphi_h)
\nonumber \\
&+\frac 1 2(\nabla\cdot \hat{{\bf u}}_h^{n+1}(\sigma^{n+1}-P_h\sigma^{n+1}),\varphi_h)
\nonumber \\
&
+(\hat{{\bf u}}_h^{n+1}\nabla e_{\sigma,h}^{n+1},\varphi_h)
+\frac{1}{2}(\nabla \cdot \hat{{\bf u}}_h^{n+1} e_{\sigma,h}^{n+1},\varphi_h)
\nonumber \\
&
+(\hat{e}_{{\bf u},h}^{n+1}\nabla\sigma^{n+1},\varphi_h)
+\frac{1}{2}(\sigma^{n+1}\nabla\cdot \hat{e}_{{\bf u},h}^{n+1},\varphi_h)
\nonumber \\
&
+((\hat{{\bf u}}^{n+1}-{\bf{R}}_h\hat{{\bf u}}^{n+1})\nabla \sigma^{n+1},\varphi_h)
\nonumber \\
&
+\frac{1}{2}(\sigma^{n+1}
\nabla \cdot(\hat{{\bf u}}^{n+1}-{\bf{R}}_h\hat{{\bf u}}^{n+1}) ,\varphi_h) \nonumber  \\
&=
(R_1^{n+1},\varphi_h),
\quad\quad\,\,
\mbox{$\forall\,\, \varphi_h\in M_h^2$.}
\end{align}
Taking
 $\varphi_h=e_{\sigma,h}^{n+1}$ in the equation above yields
\begin{align} \label{epp2}
\frac 1 {4\tau}\Big(& \|e_{\sigma,h}^{n+1}\|_{L^2}^2-\|e_{\sigma,h}^{n}\|_{L^2}^2
+\|2e_{\sigma,h}^{n+1}-e_{\sigma,h}^{n}\|_{L^2}^2
-\|2e_{\sigma,h}^{n}-e_{\sigma,h}^{n-1}\|_{L^2}^2
\Big)
\nonumber \\
\leq &
\Big|(\hat{{\bf u}}_h^{n+1}\nabla(\sigma^{n+1}-P_h\sigma ^{n+1}),e_{\sigma,h}^{n+1})\Big|
\nonumber \\
&+
\Big|\frac 1 2(\nabla\cdot \hat{{\bf u}}_h^{n+1}(\sigma^{n+1}-P_h\sigma^{n+1}),
e_{\sigma,h}^{n+1})\Big|
\nonumber \\
&+
\Big| (\hat{{\bf u}}_h^{n+1}\nabla e_{\sigma,h}^{n+1},e_{\sigma,h}^{n+1})
+\frac 1 2(\nabla \cdot \hat{{\bf u}}_h^{n+1} e_{\sigma,h}^{n+1},e_{\sigma,h}^{n+1})\Big|
\nonumber \\
&+
\Big|(\hat{e}_{{\bf u},h}^{n+1}\nabla\sigma^{n+1},e_{\sigma,h}^{n+1})
+\frac 1 2(\nabla\cdot \hat{e}_{{\bf u},h}^{n+1}\sigma^{n+1},e_{\sigma,h}^{n+1}) \Big|
\nonumber \\
&+
\Big| \Big((\hat{{\bf u}}^{n+1}-{\bf{R}}_h\hat{{\bf u}}^{n+1})
\cdot\nabla \sigma^{n+1},e_{\sigma,h}^{n+1}\Big)   \Big|
\nonumber \\
&+
\Big|\frac 1 2 \Big(\sigma^{n+1}\nabla\cdot(\hat{{\bf u}}^{n+1}
-{\bf{R}}_h\hat{{\bf u}}^{n+1}),e_{\sigma,h}^{n+1}\Big)
\Big|
\nonumber  \\
&+
\Big|(R_1^{n+1},e_{\sigma,h}^{n+1}) \Big| \nonumber \\
=&
:\sum\limits_{i = 1}^7 { I_i }.
\end{align}
By using \eqref{T1}-\eqref{T2} and \eqref{T4}, we have
\begin{align}
I_1&\leq
C \|\hat{{\bf u}}_h^{n+1}\|_{L^\infty}
\|\nabla(\sigma^{n+1}-P_h\sigma^{n+1})\|_{L^2}
\|e_{\sigma,h}^{n+1}\|_{L^2}  \nonumber \\
&\leq C \|\hat{{\bf u}}_h^{n+1}\|_{L^\infty}
h^2\|\sigma^{n+1}\|_{H^3}\|e_{\sigma,h}^{n+1}\|_{L^2}
\nonumber\\
&
\leq C\epsilon^{-1}h^4 + \epsilon\|e_{\sigma,h}^{n+1}\|_{L^2}^2,
\,\,\,\qquad\qquad\qquad\qquad\,\,\quad\qquad\quad\,\,\,
\mbox{ (use assumption \eqref{U_n}) }
\label{eph1}   \\
I_2
&\leq
C\|\nabla \hat{{\bf u}}_h^{n+1}\|_{L^\infty}
\|\sigma^{n+1}-P_h\sigma^{n+1}\|_{L^2}
\|e_{\sigma,h}^{n+1}\|_{L^2}
\nonumber\\
&\leq
C\|\nabla \hat{{\bf u}}_h^{n+1}\|_{L^\infty}
h^3\|\sigma^{n+1}\|_{H^3}
\|e_{\sigma,h}^{n+1}\|_{L^2}
\nonumber \\
&\leq
Ch^{-1}\|\hat{{\bf u}}_h^{n+1}\|_{L^\infty}
h^3
\|e_{\sigma,h}^{n+1}\|_{L^2}
\qquad\,\,\quad\,\,\qquad\qquad\quad\,
\qquad\mbox{(use \eqref{inverse})}
\nonumber \\
&\leq Ch^2\|e_{\sigma,h}^{n+1}\|_{L^2}
\qquad\,\qquad\qquad\qquad\quad\,\qquad\qquad\quad\qquad\,\,
\mbox{ (use assumption \eqref{U_n}) }
\nonumber \\
&\leq C\epsilon^{-1}h^4 + \epsilon\|e_{\sigma,h}^{n+1}\|_{L^2}^2,
\label{eph2}   \\
I_3&=0,
\label{eph3}\\
I_4
&\leq
C \|\hat{e}_{{\bf u},h}^{n+1}\|_{L^2}
\|\nabla\sigma^{n+1}\|_{L^\infty}\|e_{\sigma,h}^{n+1}\|_{L^2}
+C \|\nabla \hat{e}_{{\bf u},h}^{n+1}\|_{L^2}
\|\sigma^{n+1}\|_{L^\infty}\|e_{\sigma,h}^{n+1}\|_{L^2}
\nonumber \\
&\leq
C \left( \|\hat{e}_{{\bf u},h}^{n+1}\|_{L^2}
+\|\nabla \hat{e}_{{\bf u},h}^{n+1}\|_{L^2} \right)
\|e_{\sigma,h}^{n+1}\|_{L^2}
\nonumber \\
&\leq
C\epsilon^{-1} \left( \|\hat{e}_{{\bf u},h}^{n+1}\|_{L^2}^2
+\|\nabla \hat{e}_{{\bf u},h}^{n+1}\|_{L^2}^2  \right)
+\epsilon \|e_{\sigma,h}^{n+1}\|_{L^2}^2,
\label{eph4}\\
I_5
&\leq
C\|\hat{{\bf u}}^{n+1}-{\bf{R}}_h\hat{{\bf u}}^{n+1}\|_{L^2}
\|\nabla \sigma^{n+1}\|_{L^\infty}\|e_{\sigma,h}^{n+1}\|_{L^2}
\nonumber \\
&\leq
Ch^2\left( \|\hat{{\bf u}}^{n+1}\|_{H^2}+\|\hat{p}^{n+1}\|_{H^1}\right)
\|\nabla \sigma^{n+1}\|_{L^\infty}\|e_{\sigma,h}^{n+1}\|_{L^2}
\nonumber \\
&\leq
Ch^2\|e_{\sigma,h}^{n+1}\|_{L^2}
\nonumber \\
&\leq
C\epsilon^{-1}h^4+\epsilon\|e_{\sigma,h}^{n+1}\|_{L^2}^2,
\label{eph5}\\
I_6
& \leq
C\|\sigma^{n+1}\|_{L^\infty}\|\nabla \cdot(\hat{{\bf u}}^{n+1}
-{\bf{R}}_h\hat{{\bf u}}^{n+1})\|_{L^2}
\|e_{\sigma,h}^{n+1}\|_{L^2}
\nonumber \\
&\leq
C\|\sigma^{n+1}\|_{L^\infty}Ch^2 \left(\|\hat{{\bf u}}^{n+1}\|_{H^3}
+\|\hat{p}^{n+1}\|_{H^2}\right)
\|e_{\sigma,h}^{n+1}\|_{L^2}
\nonumber \\
&\leq
Ch^2\|e_{\sigma,h}^{n+1}\|_{L^2}
\nonumber \\
&\leq
C\epsilon^{-1}h^4+\epsilon\|e_{\sigma,h}^{n+1}\|_{L^2}^2,
\label{eph6}\\
I_7
&\leq
C\|R_1^{n+1}\|_{L^2}\|e_{\sigma,h}^{n+1}\|_{L^2}
\nonumber \\
&\leq
C\epsilon^{-1}\|R_1^{n+1}\|_{L^2}^2+\epsilon\|e_{\sigma,h}^{n+1}\|_{L^2}^2.
\label{eph7}
\end{align}
Substituting the above estimates \eqref{eph1}-\eqref{eph7} into \eqref{epp2} leads to
\begin{align}
\frac 1 {4\tau}\Big( &\|e_{\sigma,h}^{n+1}\|_{L^2}^2-\|e_{\sigma,h}^{n}\|_{L^2}^2
+\|2e_{\sigma,h}^{n+1}-e_{\sigma,h}^{n}\|_{L^2}^2
-\|2e_{\sigma,h}^{n}-e_{\sigma,h}^{n-1}\|_{L^2}^2
\Big)
\nonumber \\
\leq
&C \epsilon^{-1}\left( h^4+\|\hat{e}_{{\bf u},h}^{n+1}\|_{L^2}^2
+\|\nabla \hat{e}_{{\bf u},h}^{n+1}\|_{L^2}^2
+\|R_1^{n+1}\|_{L^2}^2      \right)
+\epsilon\|e_{\sigma,h}^{n+1}\|_{L^2}^2.
\nonumber
\end{align}
By summing up the above inequality for $n=1,2,\ldots,m$, and using \eqref{R1}, we obtain (by choosing a small $\epsilon$ )
\begin{align}
\frac 1 4\|e_{\sigma,h}^{m+1}\|_{L^2}^2
\leq
&\epsilon\tau\sum\limits_{n = 1}^{m} \|e_{\sigma,h}^{n+1}\|_{L^2}^2
+C(\tau^4+h^4) +C\tau\sum\limits_{n=0}^{m}\Big( \|e_{{\bf u},h}^{n}\|_{L^2}^2
+\|\nabla e_{{\bf u},h}^{n}\|_{L^2}^2
\Big)
\nonumber \\
&
+\frac 1 4(\|e_{\sigma,h}^1\|_{L^2}^2
+\|2e_{\sigma,h}^1-e_{\sigma,h}^0\|_{L^2}^2).
\nonumber
\end{align}
Using \eqref{F1} and  the discrete Gronwall's inequality in Lemma \ref{Gronwall's}
, we get
\begin{align}
\|e_{\sigma,h}^{m+1}\|_{L^2}^2
&\leq
C(\tau^4+ h^4) + C\tau\sum\limits_{n=0}^{m}
\left( \|e_{{\bf u},h}^{n}\|_{L^2}^2+\|\nabla e_{{\bf u},h}^{n}\|_{L^2}^2
 \right)
+C(\|e_{\sigma,h}^1\|_{L^2}^2
+\|2e_{\sigma,h}^1-e_{\sigma,h}^0\|_{L^2}^2)
\nonumber \\
&\leq
C(\tau^4+ h^4) +C\tau\sum _{n=0}^{m}
\|\nabla e_{{\bf u},h}^{n}\|_{L^2}^2
\nonumber \\
&\leq
C(\tau^4+ h^4)+C\tau\sum _{n=0}^{m}\|e_{{\bf u},h}^{n}\|_{H^1}^2,
\label{sigma_A}
\end{align}
where we use the inequality
$\|\nabla e_{{\bf u},h}^{n}\|_{L^2}\leq C\| e_{{\bf u},h}^{n}\|_{H^1}$
 to derive the  third inequality.

By the above inequality and the induction assumption \eqref{as2d}, we have
\begin{align}
&\max\limits_ {1\leq n\leq m}\|e_{\sigma,h}^{n+1}\|_{L^2}
\leq
\Big(  C(\tau^4+ h^4)+C\tau\sum _{n=0}^{m}\|e_{{\bf u},h}^{n}\|_{H^1}^2 \Big)^{\frac 1 2}
\leq
Ch(\tau^2+h^2)^{\frac 1 4}
,
\label{jie_1} \\
&\max\limits_ {1\leq n\leq m}\|e_{\sigma,h}^{n+1}\|_{L^\infty}
\leq
Ch^{-1}\max\limits_ {1\leq n\leq m}\|e_{\sigma,h}^{n+1}\|_{L^2}
\leq C(\tau^2+h^2)^{\frac 1 4}.
\label{jie 2}
\end{align}
For $\tau$ and $h$ sufficiently small, the  inequality \eqref{jie 2} implies
\begin{align}\label{rho2}
&\max\limits_ {1\leq n\leq m}\|e_{\sigma,h}^{n+1}\|_{L^\infty}
\leq \frac 1 4 \min\limits_{x\in\overline\varOmega}\sigma^0(x).
\end{align}
Then, the same argument as \eqref{sigma_n1} and \eqref{sigma_n2} shows that
\begin{align}\label{sigma_nn1}
\frac 1 2 \min_{x\in\overline\varOmega}\sigma^0(x)
\leq\sigma_h^{n+1}(x)\leq\frac 3 2 \max_{x\in\overline\varOmega}\sigma^0(x),
\quad
\mbox{$\forall\,\, x\in\varOmega,\quad n=1,\ldots,m$.}
\end{align}
\subsection{Estimates of $D_\tau e_{\sigma,h}^{n+1}$ }

Substituting
$\varphi_h=D_\tau e_{\sigma,h}^{n+1}$ in \eqref{epp1}, we have
\begin{align}\label{DD1}
\|D_\tau e_{\sigma,h}^{n+1}\|_{L^2}\leq&
\|\hat{{\bf u}}_h^{n+1}\nabla(\sigma^{n+1}-P_h\sigma^{n+1})\|_{L^2}
+\frac 1 2\|\nabla \cdot\hat{{\bf u}}_h^{n+1}(\sigma^{n+1}-P_h\sigma^{n+1})\|_{L^2}
\nonumber \\
&+\|\hat{{\bf u}}_h^{n+1}\nabla e_{\sigma,h}^{n+1}\|_{L^2}
+\frac 1 2\|\nabla \cdot\hat{{\bf u}}_h^{n+1}e_{\sigma,h}^{n+1}\|_{L^2}
\nonumber \\
&+
\|\hat{e}_{{\bf u},h}^{n+1}\nabla \sigma^{n+1}\|_{L^2}
+\frac 1 2\|\nabla \cdot\hat{e}_{{\bf u},h}^{n+1}\sigma^{n+1}\|_{L^2}
\nonumber \\
&+
\|(\hat{{\bf u}}^{n+1}-{\bf R}_h\hat{{\bf u}}^{n+1})\nabla \sigma^{n+1}\|_{L^2}
\nonumber \\
&+
\frac 1 2\|\sigma^{n+1}\nabla \cdot(\hat{{\bf u}}^{n+1}-
{\bf R}_h\hat{{\bf u}}^{n+1})\|_{L^2}
\nonumber \\
&+
\|R_1^{n+1}\|_{L^2}
\nonumber \\
=&
:\sum_{k=1}^{9} J_k,
\end{align}
where
\begin{align}
J_1+J_2\leq&
\|\hat{{\bf u}}_h^{n+1}\|_{L^\infty}\|\nabla(\sigma^{n+1}-P_h\sigma^{n+1})\|_{L^2}
+C\|\nabla \cdot \hat{{\bf u}}_h^{n+1}\|_{L^\infty}
\|\sigma^{n+1}-P_h\sigma^{n+1}\|_{L^2}
\nonumber
\\
\leq&
\|\hat{{\bf u}}_h^{n+1}\|_{L^\infty}Ch\|\sigma^{n+1}\|_{H^2}
+Ch^{-1}\|\hat{{\bf u}}_h^{n+1}\|_{L^\infty}h^2\|\sigma^{n+1}\|_{H^2}
\quad\,\,\,\,\,\,\,\quad\quad\mbox{(use \eqref{inverse})}
\nonumber \\
\leq&
Ch,
\label{D1}\\
J_3+J_4\leq&
\|\hat{{\bf u}}_h^{n+1}\|_{L^\infty}\|\nabla e_{\sigma,h}^{n+1}\|_{L^2}
+C\|\nabla \cdot \hat{{\bf u}}_h^{n+1}\|_{L^\infty}  \|e_{\sigma,h}^{n+1}\|_{L^2}
\nonumber \\
\leq&
\|\hat{{\bf u}}_h^{n+1}\|_{L^\infty}Ch^{-1}\|e_{\sigma,h}^{n+1}\|_{L^2}
+Ch^{-1}\|\hat{{\bf u}}_h^{n+1}\|_{L^\infty}\|e_{\sigma,h}^{n+1}\|_{L^2}
\nonumber \\
\leq&
\|\hat{{\bf u}}_h^{n+1}\|_{L^\infty}
Ch^{-1}h(\tau^2+h^2)^{\frac 1 4}
+Ch^{-1}\|\hat{{\bf u}}_h^{n+1}\|_{L^\infty}h(\tau^2+h^2)^{\frac 1 4}
\quad\,\,\,\,\,\,\,
\mbox{(use \eqref{jie_1})}
\nonumber \\
\leq&
C(\tau^2+h^2)^{\frac 1 4},
\qquad\qquad\quad\qquad\qquad\qquad\qquad\qquad\,\,\,
\qquad\qquad\quad\,\,\quad
\mbox{(use \eqref{U_n})}
\label{D2}\\
J_5+J_6\leq&
\|\hat{e}_{{\bf u},h}^{n+1}\|_{L^2}\|\nabla \sigma^{n+1}\|_{L^\infty}
+C\|\nabla\cdot \hat{e}_{{\bf u},h}^{n+1}\|_{L^2}
\|\sigma^{n+1}\|_{L^\infty}
\nonumber \\
\leq&
\|\hat{e}_{{\bf u},h}^{n+1}\|_{L^2}\|\nabla \sigma^{n+1}\|_{L^\infty}
+Ch^{-1}\| \hat{e}_{{\bf u},h}^{n+1}\|_{L^2}
\|\sigma^{n+1}\|_{L^\infty}
\nonumber \\
\leq&
Ch(\tau^2+h^2)^{\frac 1 4}\|\nabla \sigma^{n+1}\|_{L^\infty}
+Ch^{-1}h(\tau^2+h^2)^{\frac 1 4}\|\sigma^{n+1}\|_{L^\infty}
\quad\qquad
\mbox{(use \eqref{as2b})}
\nonumber \\
\leq&
C(\tau^2+h^2)^{\frac 1 4},
\label{D3}\\
J_7\leq&
\|\hat{{\bf u}}^{n+1}-{\bf R}_h\hat{{\bf u}}^{n+1}\|_{L^2}
\|\nabla\sigma^{n+1}\|_{L^\infty}
\nonumber \\
\leq&
Ch^2(\|\hat{{\bf u}}^{n+1}\|_{H^2}+\|\hat{p}^{n+1}\|_{H^1})
\|\nabla\sigma^{n+1}\|_{L^\infty}
\qquad\quad\quad\qquad\quad\,\,\quad\,\,\,\,\quad
\mbox{(use \eqref{T1})}
\nonumber \\
\leq&
Ch^2,
\label{D4}\\
J_8\leq&
C\|\sigma^{n+1}\|_{L^\infty}
\|\nabla\cdot(\hat{{\bf u}}^{n+1}-{\bf R}_h \hat{{\bf u}}^{n+1})\|_{L^2}
\nonumber \\
\leq&
C\|\sigma^{n+1}\|_{L^\infty}\
h(\|\hat{{\bf u}}^{n+1}\|_{H^2}+\|\hat{p}^{n+1}\|_{H^1})
\qquad\qquad\quad\quad\,\,\,\quad\qquad\,\,\,\quad
\mbox{(use \eqref{T2})}
\nonumber \\
\leq&
Ch,
\label{D5}\\
J_9\leq&
\|R_1^{n+1}\|_{L^2}
\leq
C\tau^2.
\qquad\qquad\qquad\qquad\quad\qquad\qquad\qquad\qquad\quad\,\,\,\quad
\,\,\,\mbox{(use \eqref{R1})}
\label{D6}
\end{align}
Substituting \eqref{D1}-\eqref{D6} into \eqref{DD1} yields
\begin{align}
\|D_\tau e_{\sigma,h}^{n+1}\|_{L^2}\leq&
C(\tau^2+h^2)^{\frac 1 4},
\label{D_sigma}\\
\|D_\tau e_{\sigma,h}^{n+1}\|_{L^3}\leq&
Ch^{-\frac 1 3}\|D_\tau e_{\sigma,h}^{n+1}\|_{L^2}
\qquad\quad\,\,\,\qquad\quad\qquad\quad\quad\,\,\,\quad\,\,\quad\quad\,
\mbox{(use \eqref{inverse})}
\nonumber \\
\leq&
Ch^{-\frac 1 3}(\tau^2+h^2)^{\frac 1 4}
\nonumber\\
\leq&C.
\qquad\quad\,\,\,\qquad\,\,\,
\mbox{(\,as $\tau = O(h)$ and $h$ is sufficiently small\,)}
\label{D_sigmaL3}
\end{align}

\subsection{Estimates of $e_{{\bf u},h}^{n+1}$}
The equation \eqref{eq5} and \eqref{li2} can be rewritten as
\begin{align}
&
(\sigma^{n+1}D_\tau(\sigma^{n+1}{\bf u}^{n+1}),{\bf v})
+\frac{1}{2}(\rho^{n+1}(\hat{{\bf u}}^{n+1}\cdot \nabla ) {\bf u}^{n+1}, {\bf v})
-\frac{1}{2}(\rho^{n+1}\hat{{\bf u}}^{n+1}\nabla {\bf v},{\bf u}^{n+1})\nonumber \\
&+B(({\bf u}^{n+1},p^{n+1}),({\bf v},q))=(R_2^{n+1},{\bf v}),
\qquad\quad\qquad\qquad
\forall\,\,({\bf v},q)\in \textbf{H}_0^{1}(\varOmega) \times L _0^2(\varOmega),
\label{exact_u}\\
&
(\sigma_h^{n+1}D_\tau(\sigma_h^{n+1}{\bf u}_h^{n+1}),{\bf v}_h)
+\frac{1}{2}(\rho_h^{n+1}(\hat{{\bf u}}_h^{n+1}\cdot \nabla ) {\bf u}_h^{n+1},
{\bf v}_h)
-\frac{1}{2}(\rho_h^{n+1}\hat{{\bf u}}_h^{n+1}\nabla {\bf v}_h,{\bf u}_h^{n+1})
\nonumber \\
&
+B(({\bf u}_h^{n+1},p_h^{n+1}),({\bf v}_h,q_h))=0,
\qquad\qquad\qquad\quad\qquad\,\,\,\,
\forall\,\,({\bf v}_h,q_h)\in \textbf{X}_h^{2} \times
\mathring M _h^1,
\label{sub}
\end{align}
where $R_2^{n+1}$ is truncation error.
Under the regularity assumption \eqref{smooth}, we have
\begin{align}
\tau\sum_{n=0}^{N-1}\|R_2^{n+1}\|^2_{L^2}\leq C\tau^4.
\label{TS}
\end{align}
Subtracting \eqref{sub} from \eqref{exact_u} yields
\begin{align}\label{Ee_u}
(\sigma_h^{n+1}D_{\tau}&(\sigma_h^{n+1}e_{{\bf u},h}^{n+1}),{\bf v}_h)
+B((e_{{\bf u},h}^{n+1},e_{p,h}^{n+1}),({\bf v}_h,q_h))
\nonumber \\
=
-&(\sigma_h^{n+1}D_{\tau}(e_{\sigma,h}^{n+1}{\bf u}^{n+1}),{\bf v}_h)
\nonumber \\
-&
(\sigma_h^{n+1}D_{\tau}((\sigma^{n+1}-P_h\sigma^{n+1}){\bf u}^{n+1}),{\bf v}_h)
\nonumber \\
-&
(\sigma_h^{n+1}D_{\tau}(\sigma_h^{n+1}({\bf u}^{n+1}-{\bf{R}}_h {\bf u}^{n+1}  )),{\bf v}_h)
\nonumber \\
-&
\Big[((\sigma^{n+1}-P_h\sigma^{n+1})D_{\tau}(\sigma^{n+1}{\bf u}^{n+1}),{\bf v}_h)
+(e_{\sigma,h}^{n+1}D_{\tau}(\sigma^{n+1}{\bf u}^{n+1}),{\bf v}_h)\Big]
\nonumber \\
-&
\frac{1}{2}\left[(\rho_h^{n+1}\hat{{\bf u}}_h^{n+1}\nabla e_{{\bf u},h}^{n+1},{\bf v}_h)
+(\rho_h^{n+1}\hat{{\bf u}}_h^{n+1}\nabla({\bf u}^{n+1}-{\bf{R}}_h {\bf u}^{n+1} ),{\bf v}_h)
\right]
\nonumber \\
-&
\frac{1}{2}\left[ (\rho_h^{n+1}(\hat{{\bf u}}^{n+1}-{\bf{R}}_h \hat{{\bf u}}^{n+1} )\nabla {\bf u}^{n+1},{\bf v}_h)
+(\rho_h^{n+1}\hat{e}_{{\bf u},h}^{n+1}\nabla {\bf u}^{n+1},{\bf v}_h)
\right]
\nonumber \\
-&
\frac{1}{2}\left[ ((\rho^{n+1}-P_h\rho^{n+1})\hat{{\bf u}}^{n+1}\nabla {\bf u}^{n+1},{\bf v}_h)
+(e_{\rho,h}^{n+1}\hat{{\bf u}}^{n+1}\nabla {\bf u}^{n+1},{\bf v}_h)
\right]
\nonumber \\
+&
\frac{1}{2}\left[ (\rho_h^{n+1}\hat{{\bf u}}_h^{n+1}\nabla {\bf v}_h,e_{{\bf u},h}^{n+1})
+(\rho_h^{n+1}\hat{{\bf u}}_h^{n+1}\nabla {\bf v}_h,{\bf u}^{n+1}-{\bf{R}}_h {\bf u}^{n+1} )
\right]
\nonumber \\
+&
\frac{1}{2}\left[(\rho_h^{n+1}(\hat{{\bf u}}^{n+1}-{\bf{R}}_h \hat{{\bf u}}^{n+1}) \nabla {\bf v}_h,{\bf u}^{n+1})
+(\rho_h^{n+1}\hat{e}_{{\bf u},h}^{n+1}\nabla {\bf v}_h,{\bf u}^{n+1})
\right]
\nonumber \\
+&
\frac{1}{2}\left[((\rho^{n+1}-P_h\rho^{n+1})\hat{{\bf u}}^{n+1}\nabla {\bf v}_h,{\bf u}^{n+1})
+(e_{\rho,h}^{n+1}\hat{{\bf u}}^{n+1}\nabla {\bf v}_h,{\bf u}^{n+1})
\right]
\nonumber \\
+&
(R_2^{n+1},{\bf v}_h)
\nonumber \\
=: \,\,\, &
\sum\limits_{j=1}^{11}K_j ({\bf v}_h).
\end{align}
Substituting $({\bf v}_h,q_h)=(e_{{\bf u},h}^{n+1},e_{p,h}^{n+1})$ in the above equation,
and using \eqref{T1}, \eqref{T4} and \eqref{D_sigma}, we have
\begin{align}
|K_2(e_{{\bf u},h}^{n+1})|
=&\Big|(\sigma_h^{n+1}D_{\tau}((\sigma^{n+1}-P_h\sigma^{n+1}){\bf u}^{n+1}),
e_{{\bf u},h}^{n+1})\Big|
\nonumber \\
\leq& \Big|
(\sigma_h^{n+1}D_\tau (\sigma^{n+1}-P_h\sigma^{n+1})
 {\bf u}^n,e_{{\bf u},h}^{n+1})
+\frac3 2
(\sigma_h^{n+1}(\sigma^{n+1}-P_h\sigma^{n+1})\delta_\tau{\bf u}^{n+1} ,e_{{\bf u},h}^{n+1})
\nonumber \\
&-\frac1 2 (\sigma_h^{n+1}(\sigma^{n-1}-P_h\sigma^{n-1})
\delta_\tau{\bf u}^{n},e_{{\bf u},h}^{n+1}) \Big|
\nonumber \\
\leq&
C\|\sigma_h^{n+1}\|_{L^\infty}\|{\bf u}^n\|_{L^\infty}
\|D_\tau(\sigma^{n+1}-P_h\sigma^{n+1})\|_{L^2}
\|e_{{\bf u},h}^{n+1}\|_{L^2}
\nonumber \\
&+C\|\sigma_h^{n+1}\|_{L^\infty}
\|\delta_\tau{\bf u}^{n+1}\|_{L^\infty}
\|\sigma^{n+1}-P_h\sigma^{n+1}\|_{L^2}
\|e_{{\bf u},h}^{n+1}\|_{L^2}
\nonumber \\
&+C\|\sigma_h^{n+1}\|_{L^\infty}\|\delta_\tau{\bf u}^{n}\|_{L^\infty}
\|\sigma^{n-1}-P_h\sigma^{n-1}\|_{L^2}
\|e_{{\bf u},h}^{n+1}\|_{L^2}
\nonumber \\
\leq&
Ch^2 \|D_\tau\sigma^{n+1}\|_{H^2}\|e_{{\bf u},h}^{n+1}\|_{L^2}
+Ch^2\|\sigma^{n+1}\|_{H^2}\|e_{{\bf u},h}^{n+1}\|_{L^2}
+Ch^2\|\sigma^{n-1}\|_{H^2}\|e_{{\bf u},h}^{n+1}\|_{L^2}
\nonumber \\
\leq&
Ch^4+C\|e_{{\bf u},h}^{n+1}\|_{L^2}^2,
\label{KK2}\\
|K_3(e_{{\bf u},h}^{n+1})|=&
\Big|(\sigma_h^{n+1}D_{\tau}(\sigma_h^{n+1}({\bf u}^{n+1}
-{\bf{R}}_h {\bf u}^{n+1})),e_{{\bf u},h}^{n+1})\Big|
\nonumber \\
\leq&
\Big|
(\sigma_h^{n+1}D_{\tau}e_{\sigma,h}^{n+1}
({\bf u}^n-{\bf{R}}_h {\bf u}^{n}),e_{{\bf u},h}^{n+1})\Big|
+
\Big|
(\sigma_h^{n+1}
D_{\tau}(P_h\sigma^{n+1})
({\bf u}^n-{\bf{R}}_h {\bf u}^{n}),e_{{\bf u},h}^{n+1})\Big|
\nonumber \\
&+\frac 3 2\Big|(\sigma_h^{n+1} \sigma_h^{n+1}
\delta_\tau({\bf u}^{n+1}-{\bf{R}}_h {\bf u}^{n+1})
,e_{{\bf u},h}^{n+1})\Big|
\nonumber \\
&+\frac 1 2\Big|(\sigma_h^{n+1}\sigma_h^{n-1}
\delta_\tau({\bf u}^{n}-{\bf{R}}_h {\bf u}^{n})
,e_{{\bf u},h}^{n+1})\Big|
\nonumber \\
\leq&
C\|\sigma_h^{n+1}\|_{L^\infty}\|D_\tau e_{\sigma,h}^{n+1}\|_{L^3}
\|{\bf u}^n-{\bf R}_h {\bf u}^n\|_{L^2}\|e_{{\bf u},h}^{n+1}\|_{L^6}
\nonumber \\
&
+C\|\sigma_h^{n+1}\|_{L^\infty}
\|D_{\tau}(P_h\sigma^{n+1})\|_{L^3}
\|{\bf u}^n-{\bf R}_h {\bf u}^n\|_{L^2}\|e_{{\bf u},h}^{n+1}\|_{L^6}
\nonumber \\
&+C\|\sigma_h^{n+1}\|_{L^\infty}^2
\|\delta_\tau({\bf u}^{n+1}-{\bf{R}}_h {\bf u}^{n+1})\|_{L^2}
\|e_{{\bf u},h}^{n+1}\|_{L^2}
\nonumber \\
&+C\|\sigma_h^{n+1}\|_{L^\infty}\|\sigma_h^{n-1}\|_{L^\infty}
\|\delta_\tau({\bf u}^{n}-{\bf{R}}_h {\bf u}^{n})\|_{L^2}
\|e_{{\bf u},h}^{n+1}\|_{L^2}
\nonumber \\
\leq&
Ch^2(\|{\bf u}^n\|_{H^2}+\|p^n\|_{H^1})\|\nabla e_{{\bf u},h}^{n+1}\|_{L^2}
\qquad\qquad\quad\,\,\,\qquad\qquad\qquad
\mbox{(use \eqref{D_sigmaL3})}
\nonumber \\
&+Ch^2(\|{\bf u}^n\|_{H^2}+\|p^n\|_{H^1})
\|\nabla e_{{\bf u},h}^{n+1}\|_{L^2}
\nonumber \\
&+
Ch^2(\|\delta_\tau {\bf u}^{n+1}\|_{H^2}+\|\delta_\tau p^{n+1}\|_{H^1})
\|e_{{\bf u},h}^{n+1}\|_{L^2}
\nonumber \\
&+
Ch^2(\|\delta_\tau {\bf u}^{n}\|_{H^2}+\|\delta_\tau p^{n}\|_{H^1})
\|e_{{\bf u},h}^{n+1}\|_{L^2}
\nonumber \\
\leq&
C\epsilon^{-1}h^4+C\|e_{{\bf u},h}^{n+1}\|_{L^2}^2
+\epsilon\|\nabla e_{{\bf u},h}^{n+1}\|_{L^2}^2
,
\label{KK3} \\
|K_4(e_{{\bf u},h}^{n+1})|\leq&
C\|D_{\tau}(\sigma^{n+1}{\bf u}^{n+1})\|_{L^\infty}
\|\sigma^{n+1}-P_h \sigma^{n+1}\|_{L^2}\|e_{{\bf u},h}^{n+1}\|_{L^2}
\nonumber \\
&+
C\|D_{\tau}(\sigma^{n+1}{\bf u}^{n+1})\|_{L^\infty}
\|e_{\sigma,h}^{n+1}\|_{L^2}\|e_{{\bf u},h}^{n+1}\|_{L^2}
\nonumber \\
\leq&
Ch^2\|\sigma^{n+1}\|_{H^2}
\|e_{{\bf u},h}^{n+1}\|_{L^2}
+C\|e_{\sigma,h}^{n+1}\|_{L^2}\|e_{{\bf u},h}^{n+1}\|_{L^2}
\nonumber \\
\leq&
Ch^4+C\|e_{{\bf u},h}^{n+1}\|_{L^2}^2+C\|e_{\sigma,h}^{n+1}\|_{L^2}^2,
\label{Iu4} \\
|K_5(e_{{\bf u},h}^{n+1})|\leq&
C\|\rho_h^{n+1}\|_{L^\infty}\|\hat{{\bf u}}_h^{n+1}\|_{L^\infty}
\big(\|\nabla e_{{\bf u},h}^{n+1}\|_{L^2}+\|\nabla ({\bf u}^{n+1}-{\bf R}_h {\bf u}^{n+1})\|_{L^2} \big)
\|e_{{\bf u},h}^{n+1}\|_{L^2}
\nonumber \\
\leq&
C\|\nabla e_{{\bf u},h}^{n+1}\|_{L^2}\|e_{{\bf u},h}^{n+1}\|_{L^2}
+Ch^2(\|{\bf u}^{n+1}\|_{H^3}+\|p^{n+1}\|_{H^2})\|e_{{\bf u},h}^{n+1}\|_{L^2}
\nonumber \\
\leq&
Ch^4+C\epsilon^{-1}\|e_{{\bf u},h}^{n+1}\|_{L^2}^2
+\epsilon\|\nabla e_{{\bf u},h}^{n+1}\|_{L^2}^2,
\label{Iu5}\\
|K_6(e_{{\bf u},h}^{n+1})|\leq&
C\|\rho_h^{n+1}\|_{L^\infty}\|\nabla {\bf u}^{n+1}\|_{L^\infty}
\big(\|\hat{{\bf u}}^{n+1}-{\bf R}_h \hat{{\bf u}}^{n+1}\|_{L^2}+\|\hat{e}_{{\bf u},h}^{n+1}\|_{L^2}\big)
\|e_{{\bf u},h}^{n+1}\|_{L^2}
\nonumber \\
\leq&
C \left(h^2(\|\hat{{\bf u}}^{n+1}\|_{H^2}+\|\hat{p}^{n+1}\|_{H^1})
+\|\hat{e}_{{\bf u},h}^{n+1}\|_{L^2}
\right)
\|e_{{\bf u},h}^{n+1}\|_{L^2}
\nonumber \\
\leq&
Ch^4+C\|\hat{e}_{{\bf u},h}^{n+1}\|_{L^2}^2+C\|e_{{\bf u},h}^{n+1}\|_{L^2}^2,
\label{Iu6}\\
|K_7(e_{{\bf u},h}^{n+1})|\leq&
C\|\hat{{\bf u}}^{n+1}\|_{L^\infty}\|\nabla {\bf u}^{n+1}\|_{L^\infty}
\big(\|\rho^{n+1}-P_h\rho^{n+1}\|_{L^2}+\|e_{\rho,h}^{n+1}\|_{L^2}\big)
\|e_{{\bf u},h}^{n+1}\|_{L^2}
\nonumber \\
\leq&
C(h^2\| \rho^{n+1}\|_{H^2}+\|e_{\rho,h}^{n+1}\|_{L^2})
\|e_{{\bf u},h}^{n+1}\|_{L^2}
\nonumber \\
\leq&
Ch^4+C\|e_{\rho,h}^{n+1}\|_{L^2}^2+C\|e_{{\bf u},h}^{n+1}\|_{L^2}^2,
\label{Iu7}\\
|K_{8}(e_{{\bf u},h}^{n+1})|\leq&
C\|\rho_h^{n+1}\|_{L^\infty}\|\hat{{\bf u}}_h^{n+1}\|_{L^\infty}
\|\nabla e_{{\bf u},h}^{n+1}\|_{L^2}
\big(\|e_{{\bf u},h}^{n+1}\|_{L^2}+\|{\bf u}^{n+1}-{\bf R}_h {\bf u}^{n+1}\|_{L^2}
\big)
\nonumber \\
\leq&
C
\|\nabla e_{{\bf u},h}^{n+1}\|_{L^2}
\left(\|e_{{\bf u},h}^{n+1}\|_{L^2}+Ch^2\big(\|{\bf u}^{n+1}\|_{H^2}
+\|p^{n+1}\|_{H^1}\big)\right)
\nonumber \\
\leq&
C\epsilon^{-1}(\|e_{{\bf u},h}^{n+1}\|_{L^2}^2+h^4)
+\epsilon\|\nabla e_{{\bf u},h}^{n+1}\|_{L^2}^2,
\label{Iu8}\\
|K_{9}(e_{{\bf u},h}^{n+1})|\leq&
C\|\rho_h^{n+1}\|_{L^\infty}\|{\bf u}^{n+1}\|_{L^\infty}
\big(\|\hat{{\bf u}}^{n+1}-{\bf R}_h\hat{{\bf u}}^{n+1}\|_{L^2}+\| \hat{e}_{{\bf u},h}^{n+1}\|_{L^2}\big)
\|\nabla e_{{\bf u},h}^{n+1}\|_{L^2}
\nonumber \\
\leq&
C\left(h^2\big(\|\hat{{\bf u}}^{n+1}\|_{H^2}+\|\hat{p}^{n+1}\|_{H^1}\big)
+\| \hat{e}_{{\bf u},h}^{n+1}\|_{L^2}\right)
\|\nabla e_{{\bf u},h}^{n+1}\|_{L^2}
\nonumber \\
\leq&
C\epsilon^{-1}(h^4+\| \hat{e}_{{\bf u},h}^{n+1}\|_{L^2}^2)
+\epsilon\|\nabla e_{{\bf u},h}^{n+1}\|_{L^2}^2,
\label{Iu9}\\
|K_{10}(e_{{\bf u},h}^{n+1})|\leq&
C\| \hat{{\bf u}}^{n+1}\|_{L^\infty}\|{\bf u}^{n+1}\|_{L^\infty}
\big(\|\rho^{n+1}-P_h\rho^{n+1}\|_{L^2}+\|e_{\rho,h}^{n+1}\|_{L^2}\big)
\|\nabla e_{{\bf u},h}^{n+1}\|_{L^2}
\nonumber \\
\leq&
C\big(h^2\|\rho^{n+1}\|_{H^2}+\|e_{\rho,h}^{n+1}\|_{L^2}\big)
\|\nabla e_{{\bf u},h}^{n+1}\|_{L^2}
\nonumber \\
\leq&
C\epsilon^{-1}(h^4+\|e_{\rho,h}^{n+1}\|_{L^2}^2)
+\epsilon\|\nabla e_{{\bf u},h}^{n+1}\|_{L^2}^2,
\label{Iu10}\\
|K_{11}(e_{{\bf u},h}^{n+1})|\leq&
\|R_2^{n+1}\|_{L^2}\|e_{{\bf u},h}^{n+1}\|_{L^2}
\leq
C\|e_{{\bf u},h}^{n+1}\|_{L^2}^2+C\|R_2^{n+1}\|_{L^2}^2.
\label{Iu11}
\end{align}
Moreover, the term $K_1$ can be estimated by
\begin{align}
|K_1(e_{{\bf u},h}^{n+1})|=
&\Big|(\sigma_h^{n+1}D_{\tau}(e_{\sigma,h}^{n+1}{\bf u}^{n+1}),e_{{\bf u},h}^{n+1})\Big|
\nonumber \\
\leq&
 \Big|(\sigma_h^{n+1}D_\tau e_{\sigma,h}^{n+1}\cdot {\bf u}^n,e_{{\bf u},h}^{n+1})
+\frac3 2
(\sigma_h^{n+1}e_{\sigma,h}^{n+1}\delta_\tau {\bf u}^{n+1} ,e_{{\bf u},h}^{n+1})
-\frac1 2 (\sigma_h^{n+1}e_{\sigma,h}^{n-1}\delta_\tau {\bf u}^{n},e_{{\bf u},h}^{n+1}) \Big|
\nonumber \\
\leq&
C\|\sigma_h^{n+1}\|_{L^\infty}|(D_\tau e_{\sigma,h}^{n+1}{\bf u}^n,
e_{{\bf u},h}^{n+1})|
+C\|\sigma_h^{n+1}\|_{L^\infty}\|\delta_\tau {\bf u}^{n+1} \|_{L^\infty}
\|e_{\sigma,h}^{n+1}\|_{L^2}
\|e_{{\bf u},h}^{n+1}\|_{L^2}
\nonumber \\
&+
C\|\sigma_h^{n+1}\|_{L^\infty} \|\delta_\tau{\bf u}^n \|_{L^\infty}
\|e_{\sigma,h}^{n-1}\|_{L^2}
\|e_{{\bf u},h}^{n+1}\|_{L^2}
\nonumber \\
\leq&
C|(D_\tau e_{\sigma,h}^{n+1}{\bf u}^n,e_{{\bf u},h}^{n+1})|
+C\big(\|e_{{\bf u},h}^{n+1}\|_{L^2}^2+\|e_{\sigma,h}^{n+1}\|_{L^2}^2
+\|e_{\sigma,h}^{n-1}\|_{L^2}^2\big)
\qquad\quad\,\,\,\,\,\,\,
\mbox{(use\eqref{sigma_nn1})}
\nonumber \\
\leq&
C|(D_\tau e_{\sigma,h}^{n+1}{\bf u}_h^n,e_{{\bf u},h}^{n+1})|
+C|(D_\tau e_{\sigma,h}^{n+1}({\bf u}^n-{\bf R}_h{\bf u}^n),e_{{\bf u},h}^{n+1})|
+C|(D_\tau e_{\sigma,h}^{n+1}e_{{\bf u},h}^n,
e_{{\bf u},h}^{n+1})|
\nonumber \\
&+
C\big(\|e_{{\bf u},h}^{n+1}\|_{L^2}^2+\|e_{\sigma,h}^{n+1}\|_{L^2}^2
+\|e_{\sigma,h}^{n-1}\|_{L^2}^2\big)
\nonumber \\
\leq&
C|(D_\tau e_{\sigma,h}^{n+1}{\bf u}_h^n,e_{{\bf u},h}^{n+1})|
+C\|D_\tau e_{\sigma,h}^{n+1}\|_{L^3}
\Big(\|{\bf u}^n-{\bf R}_h{\bf u}^n\|_{L^2}+\|e_{{\bf u},h}^{n}\|_{L^2}
\Big)\|e_{{\bf u},h}^{n+1}\|_{L^6}
\nonumber \\
&+
C\big(\|e_{{\bf u},h}^{n+1}\|_{L^2}^2+\|e_{\sigma,h}^{n+1}\|_{L^2}^2
+\|e_{\sigma,h}^{n-1}\|_{L^2}^2 \big)
\nonumber \\
\leq&
C|(D_\tau e_{\sigma,h}^{n+1}{\bf u}_h^n,e_{{\bf u},h}^{n+1})|
+\Big(\|{\bf u}^n-{\bf R}_h{\bf u}^n\|_{L^2}+\|e_{{\bf u},h}^{n}\|_{L^2}
\Big)\|\nabla e_{{\bf u},h}^{n+1}\|_{L^2}
\quad\quad
\mbox{(use \eqref{D_sigmaL3})}
\nonumber \\
&+
C\big(\|e_{{\bf u},h}^{n+1}\|_{L^2}^2+\|e_{\sigma,h}^{n+1}\|_{L^2}^2
+\|e_{\sigma,h}^{n-1}\|_{L^2}^2 \big)
\nonumber \\
\leq&
C\big|(D_\tau e_{\sigma,h}^{n+1}{\bf u}_h^n,e_{{\bf u},h}^{n+1})\big|
+C\epsilon^{-1}\big(h^4+\|e_{{\bf u},h}^{n}\|_{L^2}^2\big)
+\epsilon\|\nabla e_{{\bf u},h}^{n+1}\|_{L^2}^2
\nonumber \\
&+
C\big(
\|e_{{\bf u},h}^{n+1}\|_{L^2}^2
+\|e_{\sigma,h}^{n+1}\|_{L^2}^2
+\|e_{\sigma,h}^{n-1}\|_{L^2}^2\big).
\label{KKK1}
\end{align}
Now, it remains to estimate the term
$C|(D_\tau e_{\sigma,h}^{n+1}{\bf u}_h^n,e_{{\bf u},h}^{n+1})|$.
Using the definition \eqref{VV1} of $P_h$ and taking $\varphi_h =
P_h( {\bf u}_h^n \cdot e_{{\bf u},h}^{n+1}) $ in \eqref{epp1},
we get
\begin{align}\label{epp1_app}
\big|(D_\tau e_{\sigma,h}^{n+1}{\bf u}_h^n,e_{{\bf u},h}^{n+1})\big|
=&
\big|(D_\tau e_{\sigma,h}^{n+1} ,P_h( {\bf u}_h^n \cdot e_{{\bf u},h}^{n+1}))\big|
\nonumber \\
\leq&
\big|(\hat{{\bf u}}_h^{n+1}\nabla(\sigma^{n+1}-P_h\sigma ^{n+1}),P_h( {\bf u}_h^n \cdot e_{{\bf u},h}^{n+1}))\big|
\nonumber \\
&+\frac 1 2\big|(\nabla\cdot \hat{{\bf u}}_h^{n+1}(\sigma^{n+1}-P_h\sigma^{n+1}),
P_h( {\bf u}_h^n \cdot e_{{\bf u},h}^{n+1}))\big|
\nonumber \\
&
+\big|(\hat{{\bf u}}_h^{n+1}\nabla e_{\sigma,h}^{n+1},
P_h( {\bf u}_h^n \cdot e_{{\bf u},h}^{n+1}))\big|
+\frac{1}{2}\big|(\nabla \cdot \hat{{\bf u}}_h^{n+1} e_{\sigma,h}^{n+1},
P_h( {\bf u}_h^n \cdot e_{{\bf u},h}^{n+1}))\big|
\nonumber \\
&
+\big|(\hat{e}_{{\bf u},h}^{n+1}\nabla\sigma^{n+1},
P_h( {\bf u}_h^n \cdot e_{{\bf u},h}^{n+1}))\big|
+\frac{1}{2}\big|(\sigma^{n+1}\nabla\cdot \hat{e}_{{\bf u},h}^{n+1},
P_h( {\bf u}_h^n \cdot e_{{\bf u},h}^{n+1}))\big|
\nonumber \\
&
+\big|(\hat{{\bf u}}^{n+1}-{\bf{R}}_h\hat{{\bf u}}^{n+1})\nabla \sigma^{n+1},
P_h( {\bf u}_h^n \cdot e_{{\bf u},h}^{n+1}))\big|
\nonumber \\
&
+\frac{1}{2}\big|(\sigma^{n+1}
\nabla \cdot(\hat{{\bf u}}^{n+1}-{\bf{R}}_h\hat{{\bf u}}^{n+1}),
P_h( {\bf u}_h^n \cdot e_{{\bf u},h}^{n+1}))\big| \nonumber  \\
&+
\big|(R_1^{n+1},P_h( {\bf u}_h^n \cdot e_{{\bf u},h}^{n+1}))\big|
\nonumber \\
=&:\sum\limits_{i = 1}^9 { H_i }.
\end{align}
By use the \eqref{T6} and \eqref{U_n}-\eqref{U_n1}, we have
\begin{align}
H_1
\leq&
\big|(\hat{{\bf u}}_h^{n+1}\nabla(\sigma^{n+1}-P_h\sigma ^{n+1}),
P_h( {\bf u}_h^n \cdot e_{{\bf u},h}^{n+1}))\big|
\nonumber \\
\leq&
\|\hat{{\bf u}}_h^{n+1}\|_{L^\infty}
\|\nabla(\sigma^{n+1}-P_h\sigma ^{n+1})\|_{L^2}
\|P_h( {\bf u}_h^n \cdot e_{{\bf u},h}^{n+1})\|_{L^2}
\nonumber \\
\leq&
Ch^2\|\hat{{\bf u}}_h^{n+1}\|_{L^\infty}
\|\sigma^{n+1}\|_{H^3}
\|{\bf u}_h^n\|_{L^\infty}
\|e_{{\bf u},h}^{n+1}\|_{L^2}
\qquad\quad\qquad\quad\qquad\quad\qquad\qquad\quad
\mbox{(use \eqref{T4})}
\nonumber\\
\leq&
Ch^4+C\|e_{{\bf u},h}^{n+1}\|_{L^2}^2,
\label{H_1} \\
H_2
\leq&
\frac 1 2\big|(\nabla\cdot \hat{{\bf u}}_h^{n+1}(\sigma^{n+1}-P_h\sigma^{n+1}),
P_h( {\bf u}_h^n \cdot e_{{\bf u},h}^{n+1}))\big|
\nonumber \\
\leq&
C\|\nabla \cdot\hat{{\bf u}}_h^{n+1}\|_{L^\infty}
\|\sigma^{n+1}-P_h\sigma^{n+1}\|_{L^2}
\|P_h( {\bf u}_h^n \cdot e_{{\bf u},h}^{n+1})\|_{L^2}
\nonumber \\
\leq&
Ch^{-1}\|\hat{{\bf u}}_h^{n+1}\|_{L^\infty}
h^3\|\sigma^{n+1}\|_{H^3}
\|{\bf u}_h^n \|_{L^\infty}
\| e_{{\bf u},h}^{n+1}\|_{L^2}
\qquad\quad\qquad\,\,\qquad\qquad\qquad\quad
\mbox{(use \eqref{T4})}
\nonumber \\
\leq&
Ch^4+C\| e_{{\bf u},h}^{n+1}\|_{L^2}^2,
\label{H_2} \\
H_3
\leq&
\big|(\hat{{\bf u}}_h^{n+1}\nabla e_{\sigma,h}^{n+1},
P_h( {\bf u}_h^n \cdot e_{{\bf u},h}^{n+1}))\big|
\nonumber \\
=&
\Big|\left(e_{\sigma,h}^{n+1},\nabla\cdot\big(\hat{{\bf u}}_h^{n+1}
P_h( {\bf u}_h^n \cdot e_{{\bf u},h}^{n+1})\big)\right)\Big|
\nonumber \\
\leq&
C\|e_{\sigma,h}^{n+1}\|_{L^2}
\|\nabla\cdot\hat{{\bf u}}_h^{n+1}\|_{L^\infty}
\|{\bf u}_h^n\|_{L^\infty}
\| e_{{\bf u},h}^{n+1}\|_{L^2}
\nonumber \\
&+
C\|e_{\sigma,h}^{n+1}\|_{L^2}
\|\hat{{\bf u}}_h^{n+1}\|_{L^\infty}
\|\nabla\cdot{\bf u}_h^n\|_{L^\infty}
\| e_{{\bf u},h}^{n+1}\|_{L^2}
\nonumber \\
&+
C\|e_{\sigma,h}^{n+1}\|_{L^2}
\|\hat{{\bf u}}_h^{n+1}\|_{L^\infty}
\|{\bf u}_h^n\|_{L^\infty}
\| \nabla e_{{\bf u},h}^{n+1}\|_{L^2}
\nonumber \\
\leq &
C\Big(\|\nabla\hat{{\bf u}}_h^{n+1}\|_{L^\infty}
+\|\nabla{\bf u}_h^{n}\|_{L^\infty}\Big)
\|e_{\sigma,h}^{n+1}\|_{L^2}
\|e_{{\bf u},h}^{n+1}\|_{L^2}
+
C\|e_{\sigma,h}^{n+1}\|_{L^2}
\|\nabla e_{{\bf u},h}^{n+1}\|_{L^2}
\nonumber  \\
\leq&
C\Big(\|\nabla\hat{{\bf u}}_h^{n+1}\|_{L^\infty}
+\|\nabla{\bf u}_h^{n}\|_{L^\infty}\Big)
\Big(\|e_{\sigma,h}^{n+1}\|_{L^2}^2+\|e_{{\bf u},h}^{n+1}\|_{L^2}^2\Big)
+C\epsilon^{-1}\|e_{\sigma,h}^{n+1}\|_{L^2}^2
+\epsilon\|\nabla e_{{\bf u},h}^{n+1}\|_{L^2}^2,
\label{H_33}
\end{align}
Since \eqref{T3} implies
\begin{align}
\|\nabla\hat{{\bf u}}_h^{n+1}\|_{L^\infty}
+\|\nabla {\bf u}_h^{n}\|_{L^\infty}
\leq&
C\Big(\|\nabla \hat{e}_{{\bf u},h}^{n+1}\|_{L^\infty}
+\|\nabla e_{{\bf u},h}^{n}\|_{L^\infty}\Big)
\nonumber \\
&+\Big(\|\nabla ({\bf R}_h\hat{{\bf u}}^{n+1})\|_{L^\infty}
+\|\nabla ({\bf R}_h{\bf u}^{n})\|_{L^\infty}   \Big)
\nonumber \\
\leq&
Ch^{-1}\Big(\|\nabla \hat{e}_{{\bf u},h}^{n+1}\|_{L^2}
+\|\nabla e_{{\bf u},h}^{n}\|_{L^2}
\Big)+C.
\quad\qquad\,\,\,\,\,\,\,\,\,\,
\mbox{(use \eqref{inverse})}
\label{nable}
\end{align}
Substituting the above estimate into \eqref{H_33}, we have
\begin{align}
H_3\leq&
Ch^{-1}\Big(\|\nabla \hat{e}_{{\bf u},h}^{n+1}\|_{L^2}
+\|\nabla e_{{\bf u},h}^{n}\|_{L^2}
\Big)
\Big(\|e_{\sigma,h}^{n+1}\|_{L^2}^2+\|e_{{\bf u},h}^{n+1}\|_{L^2}^2\Big)
\nonumber \\
&+C\Big(\|e_{\sigma,h}^{n+1}\|_{L^2}^2+\|e_{{\bf u},h}^{n+1}\|_{L^2}^2\Big)
+
C\epsilon^{-1}\|e_{\sigma,h}^{n+1}\|_{L^2}^2
+\epsilon\|\nabla e_{{\bf u},h}^{n+1}\|_{L^2}^2.
\label{H_3}
\end{align}
Next, we have
\begin{align}
H_4
\leq&
\frac{1}{2}\big|(\nabla \cdot \hat{{\bf u}}_h^{n+1} e_{\sigma,h}^{n+1},
P_h( {\bf u}_h^n \cdot e_{{\bf u},h}^{n+1}))\big|
\nonumber \\
\leq&
C\|\nabla \cdot\hat{{\bf u}}_h^{n+1}\|_{L^\infty}
\|e_{\sigma,h}^{n+1}\|_{L^2}
\|{\bf u}_h^n \|_{L^\infty}
\|e_{{\bf u},h}^{n+1}\|_{L^2}
\nonumber \\
\leq&
C\|\nabla \hat{{\bf u}}_h^{n+1}\|_{L^\infty}
(\|e_{\sigma,h}^{n+1}\|_{L^2}^2+\|e_{{\bf u},h}^{n+1}\|_{L^2}^2)
\nonumber \\
\leq&
\Big(Ch^{-1}\big(\|\nabla \hat{e}_{{\bf u},h}^{n+1}\|_{L^2}
+\|\nabla e_{{\bf u},h}^{n}\|_{L^2}
\big)+C\Big)
\Big(\|e_{\sigma,h}^{n+1}\|_{L^2}^2+\|e_{{\bf u},h}^{n+1}\|_{L^2}^2\Big),
\quad\quad\,\,\,\,
\mbox{(use \eqref{nable})}
\label{H_4} \\
H_5
\leq&
\big|(\hat{e}_{{\bf u},h}^{n+1}\nabla\sigma^{n+1},
P_h( {\bf u}_h^n \cdot e_{{\bf u},h}^{n+1}))\big|
\nonumber \\
\leq&
C\|\hat{e}_{{\bf u},h}^{n+1}\|_{L^2}
\|\nabla\sigma^{n+1}\|_{L^\infty}
\|{\bf u}_h^n \|_{L^\infty}
\|e_{{\bf u},h}^{n+1}\|_{L^2}
\nonumber \\
\leq&
C\big(\|\hat{e}_{{\bf u},h}^{n+1}\|_{L^2}^2
+\|e_{{\bf u},h}^{n+1}\|_{L^2}^2\big),
\label{H_5} \\
H_6
\leq&
\frac{1}{2}\big|(\sigma^{n+1}\nabla\cdot \hat{e}_{{\bf u},h}^{n+1},
P_h( {\bf u}_h^n \cdot e_{{\bf u},h}^{n+1}))\big|
\nonumber \\
\leq&
C\|\sigma^{n+1}\|_{L^\infty}
\|\nabla\hat{e}_{{\bf u},h}^{n+1}\|_{L^2}
\|{\bf u}_h^n \|_{L^\infty}
\|e_{{\bf u},h}^{n+1}\|_{L^2}
\nonumber \\
\leq&
C\epsilon^{-1}\|e_{{\bf u},h}^{n+1}\|_{L^2}^2
+\epsilon\|\nabla\hat{e}_{{\bf u},h}^{n+1}\|_{L^2}^2,
\label{H_6} \\
H_7
\leq&
\big|(\hat{{\bf u}}^{n+1}-{\bf{R}}_h\hat{{\bf u}}^{n+1})\nabla \sigma^{n+1},
P_h( {\bf u}_h^n \cdot e_{{\bf u},h}^{n+1}))\big|
\nonumber \\
\leq&
C\|\hat{{\bf u}}^{n+1}-{\bf{R}}_h\hat{{\bf u}}^{n+1}\|_{L^2}
\|\nabla \sigma^{n+1}\|_{L^\infty}
\|{\bf u}_h^n \|_{L^\infty}
\|e_{{\bf u},h}^{n+1}\|_{L^2}
\nonumber \\
\leq&
Ch^4+C\|e_{{\bf u},h}^{n+1}\|_{L^2}^2,
\quad\qquad\,\,\qquad\,\,\qquad\qquad\qquad\quad
\qquad\quad\,\,\qquad\qquad\quad\,\,\,\,\quad
\mbox{(use \eqref{T1})}
\label{H_7}  \\
H_8
\leq&
\frac{1}{2}\big|(\sigma^{n+1}
\nabla \cdot(\hat{{\bf u}}^{n+1}-{\bf{R}}_h\hat{{\bf u}}^{n+1}),
P_h( {\bf u}_h^n \cdot e_{{\bf u},h}^{n+1}))\big|
\nonumber  \\
\leq&
C\|\sigma^{n+1}\|_{L^\infty}
\|\nabla(\hat{{\bf u}}^{n+1}-{\bf{R}}_h\hat{{\bf u}}^{n+1})\|_{L^2}
\|{\bf u}_h^n \|_{L^\infty}
\|e_{{\bf u},h}^{n+1}\|_{L^2}
\nonumber\\
\leq&
Ch^4+C\| e_{{\bf u},h}^{n+1}\|_{L^2}^2,
\qquad\qquad\qquad\qquad\qquad\quad\qquad\qquad
\,\,\,\qquad\qquad\quad\qquad\,
\mbox{(use \eqref{T2})}
\label{H_8} \\
H_9
\leq&
\big|(R_1^{n+1},P_h( {\bf u}_h^n \cdot e_{{\bf u},h}^{n+1}))\big|
\nonumber \\
\leq&
C\|R_1^{n+1}\|_{L^2}^2+C\| e_{{\bf u},h}^{n+1}\|_{L^2}^2.
\label{H_9}
\end{align}
Substituting \eqref{H_1}-\eqref{H_2}
and \eqref{H_3}-\eqref{H_9} into \eqref{epp1_app} yields
\begin{align}
|(D_\tau e_{\sigma,h}^{n+1}{\bf u}_h^n,e_{{\bf u},h}^{n+1})|
\leq&
Ch^4+ \epsilon\Big( \| \nabla e_{{\bf u},h}^{n+1}\|_{L^2}^2+\|\nabla\hat{e}_{{\bf u},h}^{n+1}\|_{L^2}^2\Big)
\nonumber \\
&+Ch^{-1}\Big(\|\nabla \hat{e}_{{\bf u},h}^{n+1}\|_{L^2}
+\|\nabla e_{{\bf u},h}^{n}\|_{L^2}
\Big)
\Big(\|e_{\sigma,h}^{n+1}\|_{L^2}^2+\|e_{{\bf u},h}^{n+1}\|_{L^2}^2\Big)
\nonumber\\
&+
C\epsilon^{-1}\Big(\|e_{\sigma,h}^{n+1}\|_{L^2}^2
+\|e_{{\bf u},h}^{n+1}\|_{L^2}^2
+\|\hat{e}_{{\bf u},h}^{n+1}\|_{L^2}^2
+\|R_1^{n+1}\|_{L^2}^2\Big).
\nonumber
\end{align}
Thus, we have
\begin{align}
|K_1(e_{{\bf u},h}^{n+1})|
\leq&
C|(D_\tau e_{\sigma,h}^{n+1}{\bf u}_h^n,e_{{\bf u},h}^{n+1})|
+C\epsilon^{-1}\big(h^4+\|e_{{\bf u},h}^{n}\|_{L^2}^2\big)
\nonumber\\
&+
\epsilon\|\nabla e_{{\bf u},h}^{n+1}\|_{L^2}^2
+C(\|e_{{\bf u},h}^{n+1}\|_{L^2}^2+\|e_{\sigma,h}^{n+1}\|_{L^2}^2
+\|e_{\sigma,h}^{n-1}\|_{L^2}^2)
\nonumber \\
\leq&
Ch^4+ \epsilon\Big( \| \nabla e_{{\bf u},h}^{n+1}\|_{L^2}^2+\|\nabla\hat{e}_{{\bf u},h}^{n+1}\|_{L^2}^2\Big)
\nonumber\\
&+
C\epsilon^{-1}\left(\|e_{\sigma,h}^{n+1}\|_{L^2}^2
+\|e_{\sigma,h}^{n-1}\|_{L^2}^2
+\|e_{{\bf u},h}^{n+1}\|_{L^2}^2
+\|e_{{\bf u},h}^{n}\|_{L^2}^2
+\|\hat{e}_{{\bf u},h}^{n+1}\|_{L^2}^2
+\|R_1^{n+1}\|_{L^2}^2
\right)
\nonumber \\
&+Ch^{-1}\Big(\|\nabla \hat{e}_{{\bf u},h}^{n+1}\|_{L^2}
+\|\nabla e_{{\bf u},h}^{n}\|_{L^2}
\Big)
\Big(\|e_{\sigma,h}^{n+1}\|_{L^2}^2+\|e_{{\bf u},h}^{n+1}\|_{L^2}^2\Big).
\label{K_1}
\end{align}
Then, substituting \eqref{KK2}-\eqref{Iu11} and
\eqref{K_1} into \eqref{Ee_u} yields
\begin{align}
\big|(\sigma_h^{n+1}D_\tau&(\sigma_h^{n+1}e_{{\bf u},h}^{n+1}),e_{{\bf u},h}^{n+1})\big|
+\mu\|\nabla e_{{\bf u},h}^{n+1}\|_{L^2}^2
\nonumber \\
\leq&
Ch^4+ \epsilon\left( \| \nabla e_{{\bf u},h}^{n+1}\|_{L^2}^2+\|\nabla\hat{e}_{{\bf u},h}^{n+1}\|_{L^2}^2\right)
\nonumber \\
&+Ch^{-1}\Big(\|\nabla \hat{e}_{{\bf u},h}^{n+1}\|_{L^2}
+\|\nabla e_{{\bf u},h}^{n}\|_{L^2}
\Big)
\Big(\|e_{\sigma,h}^{n+1}\|_{L^2}^2+\|e_{{\bf u},h}^{n+1}\|_{L^2}^2\Big)
\nonumber\\
&+
C\epsilon^{-1}\Big(\|e_{\sigma,h}^{n+1}\|_{L^2}^2
+\|e_{\sigma,h}^{n-1}\|_{L^2}^2
+\|e_{\rho,h}^{n+1}\|_{L^2}^2
+\|e_{{\bf u},h}^{n+1}\|_{L^2}^2
\nonumber \\
&+\|e_{{\bf u},h}^{n}\|_{L^2}^2
+\|\hat{e}_{{\bf u},h}^{n+1}\|_{L^2}^2
+\|R_1^{n+1}\|_{L^2}^2
+\|R_2^{n+1}\|_{L^2}^2\Big).
\end{align}
Summing up the above inequality  for $n=1,\ldots,m$,
and choosing sufficiently small $\epsilon$, $\tau$ and $h$,
we get
\begin{align}\label{UU00}
\max _{1\leq n\leq m}\frac 1 4\|\sigma_h^{n+1}&e_{{\bf u},h}^{ n+1}\|_{L^2}^2
+ \mu \sum_{n=1}^m\tau\|\nabla e_{{\bf u},h}^{n+1}\|_{L^2}^2
\nonumber \\
\leq&
\frac 1 4(\|\sigma_h^{1}e_{{\bf u},h}^{1}\|_{L^2}^2
+\|2\sigma_h^{1}e_{{\bf u},h}^{1}-\sigma_h^{0}e_{{\bf u},h}^{0}\|_{L^2}^2)
+C(h^4+\tau^4)+C\tau\|\nabla e_{{\bf u},h}^{0}\|_{L^2}^2
\nonumber \\
&+\tau\sum_{n=0}^{m}\left(C h^{-1}\|\nabla e_{{\bf u},h}^{n}\|_{L^2}+C \right)
\left(\|e_{\sigma,h}^{n+1}\|_{L^2}^2+\|e_{{\bf u},h}^{n+1}\|_{L^2}^2\right)
,
\end{align}
where we use \eqref{F1}-\eqref{F2}, \eqref{R1} and \eqref{TS}, and the inequality
\begin{align}
\|e_{\rho,h}^{n+1}\|_{L^2}=&\|P_h\rho^{n+1}-\rho^{n+1}+\rho^{n+1}-\rho_h^{n+1}\|_{L^2}
\nonumber \\
\leq& Ch^2+\|(\sigma^{n+1}+\sigma_h^{n+1})(\sigma^{n+1}-\sigma_h^{n+1})\|_{L^2}
\nonumber \\
\leq&Ch^2+C\|\sigma^{n+1}+\sigma^{n+1}_h\|_{L^{\infty}}\|\sigma^{n+1}-P_h\sigma^{n+1}\|_{L^2}\nonumber\\
&+C\|\sigma^{n+1}+\sigma^{n+1}_h\|_{L^{\infty}}\|P_h\sigma^{n+1}-\sigma_h^{n+1}\|_{L^2}
\nonumber \\
\leq& Ch^2+C\|e_{\sigma,h}^{n+1}\|_{L^2}.
\nonumber
\end{align}
Based on the  assumption  \eqref{F2} and \eqref{sigma_n2},
we can easily derive
\begin{align}
\|\sigma_h^{0}e_{{\bf u},h}^{0}\|_{L^2}\leq
&C(\tau^2+h^2),
\nonumber \\
\|\sigma_h^{1}e_{{\bf u},h}^{1}\|_{L^2}\leq
&C(\tau^2+h^2).
\nonumber
\end{align}
Substituting the two estimates  above into  \eqref{UU00}
and considering $\epsilon\eqref{sigma_A}+\eqref{UU00}$, we have
\begin{align}
\max_{1\leq n\leq m}\Big( \epsilon&\|e_{\sigma,h}^{n+1}\|_{L^2}^2
+\frac 1 4\|\sigma_h^{n+1}e_{{\bf u},h}^{n+1}\|_{L^2}^2 \Big)
+\mu\sum_{n=1}^m\tau\|\nabla e_{{\bf u},h}^{n+1}\|_{L^2}^2
\nonumber \\
&\leq C\epsilon\sum _{n=1}^m\tau\|\nabla e_{{\bf u},h}^{n}\|_{L^2}^2
+C(\tau^4+h^4)
\nonumber \\
&\,\,\,\,\,\,\,
+\tau\sum_{n=0}^{m}\left(C h^{-1}\|\nabla e_{{\bf u},h}^{n}\|_{L^2}+ C \right)
\left(\|e_{\sigma,h}^{n+1}\|_{L^2}^2+\|e_{{\bf u},h}^{n+1}\|_{L^2}^2\right),
\nonumber
\end{align}
where we use the fact $\tau\|\nabla e_{{\bf u},h}^0\|_{L^2}\leq C\tau h^{-1}\|e_{{\bf u},h}^0\|_{L^2}\leq C (\tau^2 +h^2)$ as $\tau = O(h)$.
By choosing $\epsilon$ small enough, the term
$C\epsilon\sum _{n=1}^m\tau\|\nabla e_{{\bf u},h}^{n}\|_{L^2}^2$
can be absorbed by the left-hand side.
Since the \eqref{as2d} implies $ \tau\sum_{n=0}^{m}\left(C h^{-1}\|\nabla e_{{\bf u},h}^{n}\|_{L^2}+C \right)\leq C$, we apply discrete Gronwall's inequality in
Lemma \ref{Gronwall's} to get
\begin{align}\label{CC2}
\max_{1\leq n\leq m}\left(\|e_{\sigma,h}^{n+1}\|_{L^2}^2
+\|e_{{\bf u},h}^{n+1}\|_{L^2}^2\right)
+\sum_{n=1}^{m}\tau\|\nabla e_{{\bf u},h}^{n+1}\|_{L^2}^2
\leq C(\tau^4+h^4),
\end{align}
where we use the upper boundedness \eqref{sigma_nn1} of $\sigma_h^{n+1}$, $n=1,2,\ldots,m$.
For sufficiently small mesh size $h$ and $\tau$, the inequality above implies
\begin{align}
\|e_{{\bf u},h}^{m+1}\|_{L^2}\leq& h(\tau^2+h^2)^{\frac 1 4},
\nonumber \\
\|e_{{\bf u},h}^{m+1}\|_{L^\infty}\leq& Ch^{-1}\|e_{{\bf u},h}^{m+1}\|_{L^2}
\leq C(\tau^2+h^2)^{\frac 1 4}\leq 1,
\nonumber \\
\tau\sum_{k=0}^{m+1}\| e_{{\bf u},h}^{k}\|_{H^1}^2
\leq&\tau\sum_{k=0}^{m+1} \|\nabla e_{{\bf u},h}^{k}\|_{L^2}^2
\leq h^2(\tau^2+h^2)^{\frac 1 2}.
\nonumber
\end{align}

Thus, the induction for \eqref{CC1a}-\eqref{CC1d} is closed.
Moreover, due to the mathematical induction on  \eqref{CC1a}-\eqref{CC1d}
is closed, it follows that \eqref{CC2} holds for $m=N-1$.
Thus, by this result and \eqref{T1}, \eqref{T4},
we complete the proof of the Theorem \ref{conclusion}.

\section{Numerical Experiments}

In this section, we present numerical examples to illustrate the convergence
and stability. The numerical results are performed by FreeFEM++.

\begin{figure}
\centering
\caption{\small Evolution of discrete energy E}
\centering
\vspace{-0.40cm}
\includegraphics[width=10.0cm,height=8.0cm]{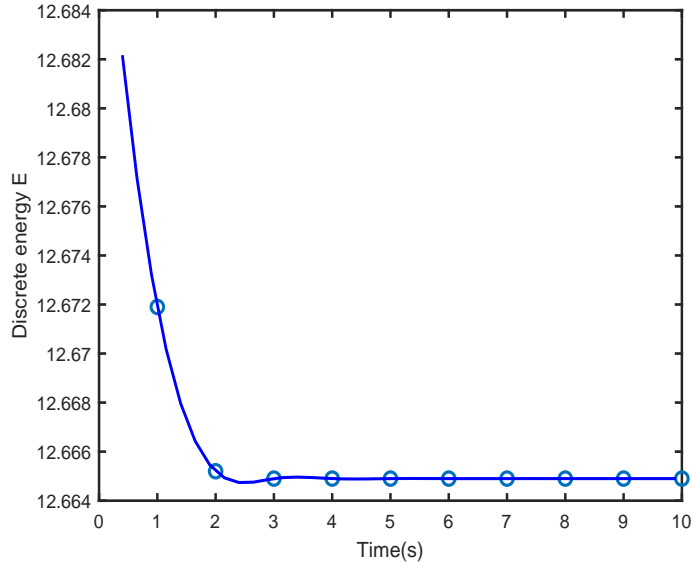}
\label{fig}
\end{figure}

\noindent
{\bf Example 4.1.}
we first consider the  follow system
\begin{align}
\partial_t \sigma + {\bf u} \cdot \nabla \sigma +\frac{1}{2} \sigma \nabla \cdot {\bf u}=g,&\label{ex11}
\\
\sigma(\sigma {\bf u})_t +\rho ({\bf u}\cdot \nabla {\bf u})
+\frac{{\bf u}}{2}\nabla\cdot (\rho {\bf u}) -\mu\Delta {\bf u} + \nabla p={\bf f}, &
\label{ex12}\\
\nabla  \cdot {\bf u}  = 0, & \label{ex13}
 \\
{\bf u} = 0,\quad \mbox{on}\quad \partial\Omega\times[0,T], &
\label{init11}\\
\rho\,(x,0) = \rho^0\quad\mbox{and}\quad {\bf u}\,(x,0) = {\bf u}^0,\,\,\hspace{0.1cm}\mbox{in} \hspace{0.1cm}\Omega, & \label{init12}
\end{align}
in a unit square domain $\varOmega=[0,1]\times[0,1]$ for $T=0.5$, where
$g$ and ${\bf f}$ are given by the following exact solutions:
\begin{equation*}
\begin{split}
\sigma(x,y,t)&=2+x(x-1)\cos(\sin(t))+y(y-1)\sin(\sin(t)) ,\\
{\bf u}(x,y,t)&=\left(\begin{array}{c}
               t^3y^2(y-1) \\
               t^3x^2(x-1)
           \end{array}
\right),\\
p(x,y,t)&=tx+y-\frac {t+1} 2.
\end{split}
\end{equation*}
The system \eqref{ex11}-\eqref{init12} is discretized by the BDF2--FEM scheme
\eqref{li1}-\eqref{li2} with the following initial conditions:
\begin{align}
&(\delta_{\tau}\sigma_h^1,\varphi_h)
+(\nabla\sigma_h^1\cdot {\bf u}_h^0,\varphi_h)
+\frac{1}{2}(\sigma_h^1\nabla \cdot {\bf u}_h^0,\varphi_h)=0,
\qquad \forall\,\,\varphi_h\in M_h^2,
\nonumber \\
&(\sigma_h^1\delta_{\tau}(\sigma_h^1 {\bf u}_h^1),{\bf v}_h)
+(\rho_h^1 {\bf u}_h^0 \nabla {\bf u}_h^1,{\bf v}_h)
\nonumber \\
&\,+\frac 1 2({\bf u}_h^1\nabla \cdot(\rho_h^1 {\bf u}_h^0),{\bf v}_h)
+B(({\bf u}_h^1,p_h^1),({\bf v}_h,q_h))=0,
\qquad \forall\,\,({\bf v}_h,q_h)\in \textbf{X}_h^{2} \times
\mathring M _h^1,
\nonumber
\end{align}
where $\rho_h^0$ and ${\bf u}_h^0$ are defined by $\rho_h^0=\Pi_h\rho^0$ and
${\bf u}_h^0={\bf{\Pi}}_h {\bf u}^0 $, respectively. To verify the convergence
order in space, we present the error results of numerical solutions with
$\tau = h$ in Table \ref{space}. We can find that the convergence rate of numerical
solutions is approximately second-order. To demonstrate the energy dissipation,
we present the numerical results of discrete energy $E$ at $T=1,2,\ldots,10$.
From the Figure \ref{fig}, we can see that the discrete energy $E$ is energy dissipation
as time goes on. These numerical results entirely consistent with theoretical analysis.

\begin{table}
  \centering
  \caption{\small Numerical errors in spatial direction ($\tau=h$)}
\begin{tabular}{c|ccccc}
\hline
$ h$&$\|\rho^N-\rho_h^N\|_{L^{2}}$&order&$\|{\bf{u}}^N-{\bf{u}}_h^N\|_{L^{2}}$ &order   \\
 \hline
$1/8$& 8.10E-03   & /       & 2.49E-05 & /       \\

$1/16$&1.98E-03   & 2.02      & 5.72E-06  &2.12    \\

$1/32$&4.85E-04  & 2.03 & 1.40E-06  & 2.02\\

$1/64$&1.20E-04  & 2.01 & 3.49E-07   & 2.00  \\
\hline
\end{tabular}\label{space}
\end{table}

\noindent
{\bf Example 4.2.}
In order to test convergence order in temporal direction, we repeat the process of Example 4.1 with  $T=1.0$. The exact solutions are provided as followings:
\begin{equation*}
\begin{split}
\sigma(x,y,t)&= 2+x(1-x)\cos(\sin(t))+y(1-y)\sin(\sin(t))  ,\\
{\bf u}(x,y,t)&=\left(\begin{array}{c}
                10x^2(x-1)^2y(y-1)(2y-1)\cos(t)  \\
               -10x(x-1)(2x-1)y^2(y-1)^2\cos(t)
           \end{array}
\right),\\
p(x,y,t)&=\sin(x)\sin(y)\sin(t).
\end{split}
\end{equation*}
To illustrate the temporal order of convergence, we solve the system \eqref{ex11}-\eqref{init12} for different step sizes $\tau$, with a fixed sufficiently small mesh size $h=\frac{1}{256}$. The errors and convergence orders of the numerical solutions are presented in Table \ref{time111}, from which we can find that our numerical scheme is second-order accurate in the temporal direction. These numerical results are  entirely consistent with the theoretical analysis in Theorem \ref{conclusion}.

\begin{table}
  \centering
  \caption{\small Numerical errors in temporal direction ($h=1/256$)}
\begin{tabular}{c|cccc}
\hline
$ \tau$&$\|\rho^N-\rho_h^N\|_{L^{2}}$&order&$\|{\bf{u}}^N-{\bf{u}}_h^N\|_{L^{2}}$ &order   \\
 \hline
$0.1$ &  4.72E-03     & /             & 7.48E-06     & /          \\

$0.1\times2^{-1}$&  1.09E-03      & 2.11       & 1.95E-06    &1.94      \\

$0.1\times2^{-2}$&  2.62E-04   &  2.06      &  5.01E-07    &1.96     \\

$0.1\times2^{-3}$&  6.40E-05    & 2.03       & 1.27E-07   &1.98     \\

$0.1\times2^{-4}$&  1.58E-05    & 2.01       & 3.21E-08   &1.99     \\
\hline
\end{tabular}\label{time111}
\end{table}

\section{Conclusions}
By introducing an equivalent form of system \eqref{ii1}--\eqref{ii4}, we propose a fully-discrete, linearized BDF2-FEM scheme for
solving the incompressible Navier--Stokes equations with density in this paper. Our numerical scheme is proved to be energy dissipation. A rigorous error estimate is presented under the assumption of sufficiently smooth of strong solutions. Finally, some numerical examples are provided to verify our theoretical analysis.

\section{Acknowledgments}

The work of Jingjing Pan and Wentao Cai was partially supported by
the Zhejiang Provincial Natural Science Foundation of China grant LY22A010019 and the National Natural Science Foundation of China grant 11901142.

\end{document}